\documentclass{cimart}

\usepackage{tikz-cd}

\newcommand{\cC}{\mathcal{C}}
\newcommand{\cE}{\mathcal{E}}
\newcommand{\cH}{\mathcal{H}}
\newcommand{\cL}{\mathcal{L}}
\newcommand{\cM}{\mathcal{M}}
\newcommand{\cS}{\mathcal{S}}
\newcommand{\cT}{\mathcal{T}}
\newcommand{\cU}{\mathcal{U}}
\newcommand{\C}{\mathbb{C}}
\newcommand{\N}{\mathbb{N}}
\newcommand{\R}{\mathbb{R}}
\newcommand{\Z}{\mathbb{Z}}

\DeclareMathOperator{\Spec}{Spec}
\DeclareMathOperator{\Tr}{Tr}

\title{On quasi-isospectrality of potentials and Riemannian mani\-folds}

\authors{Clara L. Aldana and Camilo Pérez}

\authorinfo[C. L. Aldana]{Independent Researcher, Regensburg, Germany}{clara.aldana@posteo.net}

\authorinfo[C. Pérez]{Universidad de los Andes, Bogotá, Colombia}{ca.perezt@uniandes.edu.co}

\abstract{%
    In this article, we study quasi-isospectral operators as a generalization of isospectral operators. The paper contains both expository material and original results. We begin by reviewing known results on isospectral potentials on compact manifolds and finite intervals, and then introduce the notion of quasi-isospectrality. We next investigate the BMT method as a systematic approach to constructing quasi-isospectral Sturm–Liouville operators on a finite interval, and apply it to several boundary value problems.
    Our main result shows that any two quasi-isospectral closed manifolds of odd dimension are, in fact, isospectral. In addition, we extend classical compactness results for isospectral potentials on low-dimensional manifolds to the quasi-isospectral setting via heat trace asymptotics.
    }



\VOLUME{34}
\YEAR{2026}
\ISSUE{2}
\NUMBER{4}
\DOI{https://doi.org/10.46298/cm.15976}
\licence{arXiv.org - Non-exclusive}

\editinfo{July 2, 2025}{December 15, 2025}{Tiago Macedo, Jaqueline Mesquita, Rafael Potrie, and Mariel Saez}

\acknowledgments{
\noindent The authors acknowledge the support of the ``Fundaci\'on para la promoci\'on de la investigaci\'on y la tecnolog\'{\i}a'' of the Banco de la Rep\'ublica de Colombia through project number~4.498, contract number~202016. 
We are grateful to the two anonymous referees for pointing out an inaccuracy in the statement of Theorem~\ref{thm:qiRmnb-hi}, which led to an improvement of the results presented in this paper, and for their constructive comments that helped enhance its overall quality.\\
This article is the result of work that began in the second semester of 2020, during the Covid-19 pandemic. We would like to thank all those around us whose support helped us through this challenging period. Finally, we are grateful to the editor Mariel Saez, for suggesting that we consider this volume for submission.

}

\begin{document}

\section{Introduction}

The isospectral problem is a classical problem in mathematics. It poses the following question: if the spectra of two operators are equal, are the operators the same up to unitary transformations? In the geometric setting, this question was made famous by Mark Kac's formulation of it: \lq\lq Can we hear the shape of a drum?\rq\rq \  It refers to the relation between the frequencies that we hear of a string and the spectrum of the Laplace Operator defined on functions on the string.  Isospectral manifolds have been widely studied since then, with many positive and negative results.

Two closed Riemannian manifolds are called isospectral if the spectra of the corresponding Laplace operators are the same, including multiplicities. If the manifold is compact with a boundary, a boundary condition should be specified in order to obtain a closed extension of the Laplacian.

The setting is a compact Riemannian manifold $(M,g)$ of dimension $n$, possibly with boundary, and a closed extension $\Delta$ of the associated Laplacian. In coordinates, we write $g = (g_{ij})_{i,j =1, \dots, n}$ and
\[ \Delta_g = - \frac{1}{\sqrt{\det(g)}} \sum_{i,j=1} ^n \partial_i g^{ij} \sqrt{\det(g)} \partial_j. \]
Notice that with this sign convention, the Laplacian is a positive operator. If we need to emphasize the dependence on the metric, we keep the subindex $g$, otherwise we may drop it. Let us consider a function $V\in \cC^{\infty}(M)$, or possibly $V\in \cL^p(M)$ for some $p\geq1$. This function is called the potential, and
$$\Delta +V$$
is the Schr\"odinger operator associated with it. Two potentials are called isospectral if the spectra of the corresponding Schr\"odinger operators are the same, including multiplicities. Isospectral potentials have been extensively studied.

This article contains both expository material and original results. At the beginning of each section, we indicate whether it is expository or original. Some of the expository sections also include observations and comments by the authors. In Section~\ref{sec:knownresults} we present some relevant rigidity and compactness results for isospectral potentials. In Section~\ref{sec:StLiD}, we follow \cite{PoeschTrub} to explain the well-understood structure of the set of isospectral potentials for the Sturm-Liouville problem on a closed interval with Dirichlet boundary conditions.

The spectral theory of Schr\"odinger operators is fundamental in quantum mechanics, where the spectrum corresponds to the possible energy levels of a quantum system. In Mathematical Physics, the interplay between the potential and the spectral data governs a wide range of inverse problems, from determining atomic structures to understanding wave propagation in inhomogeneous media. Precise spectral information also underpins semiclassical analysis and the study of quantum chaos, linking spectral data to the underlying classical dynamics.

Before we continue, let us introduce the heat semigroup
$e^{-t\Delta_g}$, for $t>0$. It is used to construct solutions to the heat equation. If the manifold $M$  has a boundary, and we consider the initial value problem with Dirichlet boundary conditions
 \[(\partial_t +\Delta_g) u =0, u(z,0)=f(z), \left.u(z,t)\right\vert_{\partial M} =0,\]
For each $t>0$, the heat operator is a compact, smoothing operator, and it has an integral kernel $K_g(z,z',t)$. In terms of this, the solution to the heat equation can be written as
$$u(z,t) = (e^{-t\Delta_g}f)(z)=\int_{M} K_{g}(z,z',t) f(z')\ d\mu(z')$$
when $f \in \cL^2(M)$.
Moreover, the heat operator is trace class for every  $t>0$, if we denote the spectrum of the Laplace operator as  $\Spec(\Delta_g) =(\lambda_j )_{j \in \N}$, we have that
\begin{equation}
\Tr(e^{-t\Delta_g}) = \sum_{j=0}^{\infty} e^{-t\lambda_j} = \int_M K_g(z,z,t)d\mu(z).
\label{eq:def-ht}
\end{equation}
As we will see in this article, the heat trace plays an important role in spectral and inverse spectral theory. If two operators are isospectral, they have the same heat trace. Furthermore, the spectrum of $\Delta_g$ and the heat trace, viewed as a function of $t\in ]0,\infty[$, completely determine each other \cite[p.~132]{craioveanu2013old}, \cite[p.~120]{lablee2015spectral}.

As the parameter $t$ tends to zero, the heat operator goes to the identity, $e^{-t\Delta_g} \to Id$, and the heat trace diverges.
After Minakshisundaram and Pleijel formulas, see \cite{Minakshisundaram1949SomePO}, \cite[p.~215]{berger1971spectre},  \cite[p.~99]{rosenberg1997laplacian}, there is an asymptotic expansion of the heat kernel that implies an asymptotic expansion of the heat trace.
If $M$ is closed, as  $t\to 0^+$
\begin{eqnarray*}
K_g (z,z',t) & \sim & t^{-n/2} \sum_{j=0}^\infty u_j(z,z') t^j,\\
\Tr(e^{-t\Delta_{g}}) &\sim & t^{-n/2}\sum_{j=0}^\infty a_{j} t^{j}.
\end{eqnarray*}

If on the contrary $M$ has a boundary, $\partial M \ne \emptyset$, we have that
\[
\Tr(e^{-t\Delta_{g}}) \sim t^{-n/2}\sum_{j=0}^\infty a_{j} t^{j/2}.
\]

The coefficients in the expansion are  called the heat invariants. The heat invariants of isospectral metrics are the same. If $M$ is closed,
$$a_{j} =\int_{M} u_{j}(z,z) d\mu(z).$$
If $M$ has a boundary, the heat invariants are specific sums of integrals over the interior of $M$ and integrals over the boundary. The latter are of the form:
$$a_k = \int_{\partial M} v_k(w,w) d\nu(w).$$

\begin{example} For a $3$-dimensional closed manifold, the first three heat invariants are the following:
$$a_0 = \frac{\operatorname{Vol}(M)}{(4\pi)^{3/2}}.$$
If $R$ denotes the scalar curvature,
$$a_1 = \frac{1}{(4\pi)^{3/2}} \cdot \frac{1}{6} \int_M R \, \mathrm{dV}_g.$$
Denoting the Ricci curvature by $\operatorname{Ric}$ and the Riemann tensor by $\operatorname{Rm}$
$$a_2 = \frac{1}{(4\pi)^{3/2}} \cdot \frac{1}{360} \int_M \left( 5R^2 - 2|\operatorname{Ric}|^2 + |\operatorname{Rm}|^2 \right) \mathrm{dV}_g.$$
In dimension $2$, the first two heat invariants are
$$a_{0}=\text{Area}(M,g), \quad \quad a_1= \frac{\chi(M)}{6}.$$
Hence, the spectrum determines the Area and the Euler characteristic.
\end{example}

Our main motivation was to study generalisations of the concept of isospectrality and to see which of the properties of isospectral sets of potentials still hold for these weaker notions. In this paper, we consider the concept of quasi-isospectral potentials as it was introduced by Bilotta et al in \cite{BILOTTA20171050}. However, here we introduce the notion of quasi-isospectrality in a slightly more general setting.
\begin{definition}
Let $T_1$, $T_2$ be two closed self-adjoint linear operators defined on Hilbert spaces  $\cH_1$, $\cH_2$, respectively. Assume that each of them has a compact resolvent and a discrete spectrum. We say that $T_1$ and $T_2$ are {\it quasi-isospectral} if their spectra are the same, with the exception of one of the eigenvalues. The two differing eigenvalues must lie in a prescribed interval of fixed length that contains no other eigenvalues of either operator.
\end{definition}

Refinements of this definition are given depending on the concrete setting. For example, if we consider the Sturm-Liouville eigenvalue problem with Dirichlet boundary conditions,
\begin{equation}
\begin{cases}
- y '' (x) + q(x) y(x) = \lambda y(x), & \quad x\in ]0,1[ \\
y(0) = 0 = y(1), &
\end{cases} \label{SLp-Dirichlet}
\end{equation}
with $q$ a continuous function on $[0,1]$, and if we denote the corresponding operator by $\cT_q$, we can define quasi-isospectral potentials as follows:

\begin{definition}
Two potentials $q_1$ and $q_2$ are called quasi-isospectral if the corresponding operators $\cT_{q_1}$ and $\cT_{q_2}$ are quasi-isospectral.
\end{definition}

In Section~\ref{sec:Quasi-Iso} we explain the construction of quasi-isospectral potentials on finite intervals. A natural question that arises is whether there are potentials that are isospectral, but the boundary conditions are different. In Section~\ref{sec:Examples}, we exhibit two examples in which this happens. In Section~\ref{s:cimap}, we briefly review some compactness results for Riemannian manifolds and for Schr\"odinger operators on closed manifolds and planar domains.

Finally, in sections~\ref{s:qiRm} and~\ref{sec:Results}, which contain most of the original results of this paper, we apply the concept of quasi-isospectrality to compact Riemannian manifolds and to potentials on compact manifolds, respectively. The main result of this part is Theorem~\ref{thm:qiRmnb-hi}, where we establish that quasi-isospectral closed manifolds of odd dimension are, in fact, isospectral. We also generalize several compactness results that hold for isospectral potentials on closed manifolds to the quasi-isospectral setting. We obtain these results through an analysis of the heat trace and the heat invariants. Finally, for quasi-isospectral Sturm–Liouville operators on a closed interval with Dirichlet boundary conditions, we show that the mean value of the potentials over the interval must be the same.

A brief note on the notation: in some parts of this paper, $g$ denotes a Riemannian metric on a manifold, while in others, the same letter $g$ denotes an eigenfunction of a boundary value problem on an interval. Although this may seem confusing at first glance, the context in each section makes the intended meaning of $g$ clear. The same holds for the letter $h$.

\section{Review of known rigidity results}
\label{sec:knownresults}
In this section, we present some classical results about rigidity of Schr\"odinger operators that are relevant for this paper, with some additional observations at the end of Subsection~\ref{ss:srsocm}.

\subsection{Sturm-Liouville operators on a finite interval}

One key point in inverse spectral theory is to know whether there exist potentials that are indeed isospectral. The theorems presented in this section are related to the opposite case, that is, they give conditions under which the spectrum of the Schr\"odinger operator determines the corresponding potential.
The presentation given here follows \cite{Freiling2001InverseSP}. Let us start with a simple example.

\begin{example}
	Let us consider the Laplace operator $\Delta_{\mathcal{I}}u = -u''$ on a {\emph string}, $\mathcal{I} =[0,1]$.\label{ex:eigenspace}
\begin{enumerate}
\item  The Dirichlet eigenvalue problem associated with the Laplacian on $\mathcal{I}$
$$\Delta_{\mathcal{I}} u = \lambda u, \quad \quad u(0)=0=u(1), \quad \quad u\neq 0,$$
has solutions
$$u_n(t)=\sin(n\pi t), \text{ with } \lambda_n = n^2 \pi^2  \text{ and } n\geq 1.$$
 The spectrum of the Dirichlet Laplacian is
 $\Spec(\Delta_{\mathcal{I}_D})=(n^2 \pi^2)_{n\geq 1}$.
\item  The Neumann eigenvalue problem associated to $\Delta_{\mathcal{I}}$
$$\Delta_{\mathcal{I}} u = \lambda u, \quad \quad u'(0)=0=u'(1), \quad \quad u\neq 0,$$
has solutions
$$u_n(t)=\cos(n\pi t), \text{ with } \lambda_n = n^2 \pi^2 \text{ and } n\geq 0.$$
 The spectrum of the Neumann Laplacian is
 $\Spec(\Delta_{\mathcal{I}_N})=(n^2 \pi^2)_{n\geq 0}$.
\end{enumerate}
\end{example}

The following theorem by Ambarzumian states that for the Neumann Laplacian on the unit interval, there are no isospectral potentials.

\begin{theorem}[Ambarzumian's theorem \cite{Ambarzumian}]
If the eigenvalues of the Neumann boundary value problem (BVP)
$$\begin{cases}
- y '' (x) + q(x) y(x) = \lambda y(x), & \quad x\in ]0,1[ \\
y'(0) = 0 = y'(1) &
\end{cases}
$$
are $\lambda_n = \pi^2 n^2$, $n\geq0$, then $q(x)=0$ a.e. on $]0,1[$.
\end{theorem}

Ambarzumian's theorem is an exception since, in general, the set of potentials isospectral to a given one is not empty. However, imposing additional restrictions, there are  positive results for certain instances of Robin boundary conditions, as the next two theorems show.

\begin{theorem}[Marchenko\cite{Marchenko}-Levinson\cite{levinson}]
Let $(\lambda_n)_{n\geq0}$ be the eigenvalues of the BVP
$$\begin{cases}
- y '' (x) + q(x) y(x) = \lambda y(x), & \quad x\in ]0,1[ \\
y'(0) -h y(0) = 0, \quad y'(1) + H y(1) =0 &
\end{cases}
$$
with $h,H \in \R$. Let $\alpha_n = \| \varphi(\cdot, \lambda_n) \|_{\cL^2[0,1]}$
where $\varphi(x, \lambda) $ is a solution of
$$- y '' (x) + q(x) y(x) = \lambda y(x)$$
that satisfies $\varphi(0, \lambda) =1$, $\varphi'(0, \lambda)=h$. Then $(\lambda_n,\alpha_n)_{n\geq 0}$ determine the potential $q(x)$ a.e. and the constants $h, H$. If, in addition, the potential has the symmetry
$$q(1-x) = q(x), \quad \quad \forall x\in [0,1],$$ so $H=h$, then we have that the spectrum $(\lambda_n)_{n\geq 0}$ determines $q(x)$ and the value $h$.
  \end{theorem}

The Marchenko-Levinson theorem shows that a potential is uniquely determined a.e. by the spectrum of one closed extension of the operator (via boundary conditions) if additional data, such as normalization constants, are provided. However, a key question in inverse spectral theory is whether uniqueness can be achieved using only spectral data. In this context, Borg's classical result, given below, shows that the spectra of two closed extensions corresponding to two distinct boundary conditions are sufficient to determine the potential uniquely.

\begin{theorem}[Borg \cite{Borg}]
Let $(\lambda_n)_{n\geq0}$ be the eigenvalues of the BVP
$$\begin{cases}
- y '' (x) + q(x) y(x) = \lambda y(x), & \quad x\in ]0,1[ \\
y'(0) -h y(0) = 0, \quad y'(1) + H y(1) =0 &
\end{cases}
$$
with $h,H \in \R$. Let $(\mu_n)_{n\geq 0}$ be the eigenvalues of the BVP
$$\begin{cases}
- y '' (x) + q(x) y(x) = \lambda y(x), & \quad x\in ]0,1[ \\
y'(0) -h y(0) = 0,  \quad y(1) =0, &
\end{cases}
$$
then $(\lambda_n,\mu_n)_{n\geq 0}$ determines the potential $q(x)$ a.e. and the constant $h$.
\end{theorem}

\subsection{Spectral rigidity of Schr\"odinger operators on compact manifolds.}
\label{ss:srsocm}
In \cite{Davies2010}, E. B. Davies proved an extension of Ambarzumyan's theorem to the setting of a compact metric measure space $X$ for a non-negative self-adjoint operator $H_0$ acting on $\cL^2(X,dx)$. Davies imposed four conditions on the heat semigroup associated with the operator $H_0$. These conditions are satisfied by the heat operator of the Laplacian on a closed manifold, or the one of the Laplacian with Neumann boundary conditions on a compact manifold with a boundary that satisfies the Lipschitz condition. In the case when $X$ is a compact connected Riemannian manifold and $H_0$ is a closed extension of the Laplacian, which we denote by $\Delta$, the conditions are:

\begin{itemize}
    \item The operator $e^{-t \Delta}$ has a non-negative kernel $K(t,x,y)$ which is continuous on $(0,\infty)\times X\times X$.
    \item There exist constants $c,d>0$ such that $0\leq K(t,x,x)\leq ct^{-d/2}$, for all $t\in]0,1[$.
    \item There exist a constant $a>0$ such that $t^{d/2}K(t,x,x)\to a$ when $t\to0$. If $\partial X\neq\emptyset$, the previous limit holds for all $x\notin \partial X$.
    \item The smallest eigenvalue is $\lambda_1=0$ with multiplicity one. The corresponding normalized eigenfunction is constant with value $\phi_1=\text{Vol}(X)^{-1/2}$.
\end{itemize}

Now, let $H_V:=\Delta+V$, where $V$ is a bounded real-valued potential, or more generally $V\in L^{\infty}_{loc}(X)$. It is well known that $\Spec(H_V)$ is discrete and has
$\inf_{x\in X}V(x)$ as lower bound, see \cite[p.101]{lablee2015spectral}.

If we denote $\Spec(\Delta)=(\lambda_k)_{k\geq1}$ and $\Spec(H_V)=(\mu_k)_{k\geq1}$, from the min-max principle, the eigenvalues are related as follows:

\begin{equation}
  \lambda_k-\|V\|_{\infty}\leq \mu_k\leq \lambda_k+\|V\|_{\infty} ,
  \label{compara}
\end{equation}
for all $k\geq1$. A straightforward computation using equation (\ref{compara}) gives the comparison of the heat traces:
$$e^{-t \|V\|_{\infty}}\;\text{Tr}\left(e^{-t \Delta}\right)\leq\text{Tr}\left(e^{- t H_V}\right)\leq e^{t \|V\|_{\infty}}\;\text{Tr}\left(e^{-t \Delta}\right).$$

In \cite[Section 4]{Davies2010}, Davies states and proves Theorems 5, 6, and 7, which we summarize in the following theorem:

\begin{theorem}[Davies]
With the notation introduced above, if $\mu_1\geq0$ and any of the following conditions is satisfied
\begin{enumerate}
\item $\int_{X} V dx \leq 0,$
\item $\limsup_{t\to0} t^{d/2-1}\sum_{n=1}^{\infty}\left(e^{-t \lambda_n}-e^{-t \mu_n}\right) \leq 0,$ or
\item $\limsup_{n\to0}(\mu_n-\lambda_n)\leq0$,
\end{enumerate}
we have that $V=0$.
\end{theorem}

Notice that Davies assumes that $ \lambda_{1}=0 $ with multiplicity one. He also proves that the condition $ \int_{x\in X}V(x)\le 0 $ is implied by both $(2)$ and $(3)$. Moreover, $(1)$ implies that the first eigenvalue of $H_{V} $ must be non-positive since \[ \mu_{1}\le \text{Vol}^{-1}(X)\int_{x\in X}V(x)\le 0, \]
where if $ \phi_{1}=\text{Vol}^{-1/2}(X) $ is a $ \lambda_{1}$-eigenfunction, due to the characterization of the first eigenvalue given by the min-max principle.

We now would like to mention the implications of Davies’s results for potentials that are isospectral or quasi-isospectral to the Laplacian. Assume that $\Spec(\Delta+V)=\Spec(\Delta)$, including multiplicities. Then, $\mu_1=\lambda_1=0$ and condition~(3) holds with equa\-lity, hence $V\equiv 0$. Now, if the operators are instead quasi-isospectral, so $\Spec(\Delta+V)$ and $\Spec(\Delta)$ differ only at finitely many eigenvalues, and if $\mu_1\ge 0$, then $\limsup_{n\to\infty}(\mu_n-\lambda_n)=0$, and Davies’ criterion again forces $V\equiv 0$.

\section{No rigidity of Sturm-Liouville with Dirichlet boundary conditions}
\label{sec:StLiD}
In this section, we summarize known facts concerning classical inverse spectral theory of Sturm-Liouville operators with Dirichlet boundary conditions on the interval $[0,1]$, following
\cite{PoeschTrub, Freiling2001InverseSP}. We consider sets of isospectral potentials in {$\cL^2[0,1]$}. These sets are known to be non-trivial, and their structure is well understood.

Let us recall the boundary value problem from equation (\ref{SLp-Dirichlet}):
$$\begin{cases}
- y '' (x) + q(x) y(x) = \lambda y(x), & \quad x\in ]0,1[\\
y(0) = 0 = y(1), &
\end{cases}  $$
but now we consider $q\in \cL^2[0,1]$, and let $\cT_q$ be the corresponding operator. Following the shooting method, let $y_2(x,\mu)$ be a solution of the problem
$$\begin{cases}
- y '' (x) + q(x) y(x) = \mu y(x), & \quad x\in ]0,1[ \\
y(0) = 0, y'(0) = 1. & \\
\end{cases}
$$
The zeros $\lambda$ of the map $\mu\mapsto y_2(1,\mu)$ are eigenvalues of the Dirichlet problem (\ref{SLp-Dirichlet}) and $y_2(x,\lambda)$ is a $\lambda$- eigenfunction. In addition,
\begin{itemize}
    \item The $\lambda$-eigenspace is one dimensional. Any other $\lambda$-eigenfunction is scalar multiple of $y_2(x,\lambda)$. In particular $\varphi(x)=\varphi'(0)y_2(x,\lambda)$.
    \item {$\lambda$ has algebraic multiplicity one}. The Dirichlet spectrum is the zero set of the function $y_2(1,\cdot)$. On the other hand, if $\lambda$ is a Dirichlet eigenvalue, then
    $$\left[\dfrac{\partial}{\partial \lambda}y_2(1,\lambda)\right]\left[\dfrac{\partial}{\partial x}y_2(x,\lambda)\bigg|_{x=1}\right]=\|y_2(\cdot,\lambda)\|^2$$
    which implies that the roots of $y_2$ are simple, see \cite[Thm.2 Ch.2]{PoeschTrub}.
\end{itemize}

Notice that each eigenvalue $\lambda_n(q)$ of $\cT_q$, has a unique eigenfunction $z_n(x,q)$ such that $\|z_n\|=1$ and $\frac{d}{dx}z_n(0,q)>0$. Moreover, as a function of the potential, $\lambda_n$ is a real analytic function on $\cL^2[0,1]$, and $\frac{\partial\ }{\partial q(t)}\lambda_n=z_n^2(t,q)$, see
\cite[Ch.2.]{PoeschTrub}.

The eigenvalues have the following asymptotic expansion as $n\to \infty$,
$$\lambda_n=n^2\pi^2+ \int_0^1 q(t)\;dt - \langle \cos(2\pi n x) , q\rangle + O\left(\frac1n \right).$$
From which we obtain the very useful expression for eigenvalues,

\begin{equation}\label{explicit}
\lambda_n(q)
= n^2\pi^2 + C(q) + b_n(q),\qquad
C(q):=\int_0^1 q(t)\,dt,\quad
b_n(q)=\lambda_n(q)-n^2\pi^2-C(q),
\end{equation}
where the sequence $b(q):=(b_n(q))_{n\ge1}$ belongs to $\ell^2$, see \cite[Thm.4, Ch.2]{PoeschTrub}. Notice that the term $C(q)$ is independent of $n$.

Next, we consider the maps:
$$\begin{tikzcd}
\cL^2[0,1] \ar{r}{\Lambda} & \cS \subset \R^{\N} \ar{r}{\Omega} & \mathbb{R}\times\ell^2\\\vspace{-5mm}
 q  \arrow[r, maps to] & (\lambda_n(q))_{n\geq0} \arrow[r, maps to] & (C(q),b(q)).
\end{tikzcd}$$
Here $\Lambda(q) = \Spec(\cT_q)$, and $\cS$ is the set of sequences of eigenvalues that have an asymptotic expansion as described above. The map $\Omega$ is injective, and the map $\Lambda$ is real analytic from $\cL^2[0,1]$ into $S$. In particular, it is well known that, for every $v\in \cL^{2}[0,1]$, the maps
$$\varepsilon\mapsto \lambda_n(q+\varepsilon v), \ \text{ and } \ \varepsilon\mapsto z_n(\cdot,q+\varepsilon v)$$
are analytic near $\varepsilon=0$.
Write $\lambda_n(\varepsilon)=\lambda_n(q+\varepsilon v)$ and $z_n(\varepsilon)=z_n(\cdot,q+\varepsilon v)$
with the normalization $\|z_n(\varepsilon)\|_{2}=1$. If $ \cM_{v} $ is the operator multiplication by $ v $, differentiating the eigenvalue equation
\[ (\cT_q+\varepsilon \cM_v)z_n(\varepsilon)=\lambda_n(\varepsilon)\,z_n(\varepsilon) \]
with respect to $\varepsilon$, evaluating at $ \varepsilon =0$,  and taking the inner product with $ z_{n}(0)=z_{n}(\cdot,q)$, we obtain:
\[ \langle z_n,\,\cT_q \dot{z}_n(0)\rangle+\langle z_n,\,\cM_v z_n\rangle
=\dot{\lambda}_n(0)\langle z_n,z_n\rangle+\lambda_n(q)\langle z_n,\dot{z}_n(0)\rangle,\]
where the dot indicates the derivative with respect to $\varepsilon$.
Since $\cT_q$ is self-adjoint and $\cT_q z_n=\lambda_n(q) z_n$, yields the formula:
\[D \lambda _{n}(q)(v) := \left.\frac{d\ }{d\varepsilon} \lambda_{n}(q + \varepsilon v) \right|_{\varepsilon = 0} = \langle z_{n}^{2},v \rangle =\int_0^1 v(x)\,z_n(x,q)^2\,dx.\]
Since $DC(q)[v] := \left.\frac{d\ }{d\varepsilon} C(q + \varepsilon v) \right|_{\varepsilon = 0} = \langle 1,v\rangle$, we have
\[
Db_n(q)( v ) \;=\; D\lambda_n(q)(v)-DC(q)(v) \;=\; \langle z_n^2-1,\,v\rangle.
\]
Therefore,
\[
D\Lambda_q(v)=\big(C(v),\,(\,\langle z_n^2-1,\,v\rangle\,)_{n\ge1}\big),
\]
see \cite[Thm.1 Ch.3]{PoeschTrub}. The analysis that follows consists of dividing into the cases of even and odd potentials. Let
$$\cU:=\{q\in \cL^2[0,1]\;|\; q(1-x)=-q(x);\; x\in[0,1]\}$$
be the subspace of all odd functions and
$$\cE:=\{q\in \cL^2[0,1]\;|\; q(1-x)=q(x);\; x\in[0,1]\}$$
be the subspace of all even functions. $\cE$ is the orthogonal complement of $\cU$ in $\cL^2[0,1]$, so it decomposes as $\cL^2[0,1]=\cE \oplus \cU$.

An easy computation shows that for $q=0$, $z_n^2(x,q)=2\sin^2{(\pi n x)}$, and it is not difficult to prove that the kernel of $D\Lambda_0$ is the space of all functions whose Fourier cosine  coefficients vanish. This implies that $D\Lambda_0$ is a linear isomorphism when restricted to the space of even functions $\cE$, \cite[p.53]{PoeschTrub}.
Therefore, $\Lambda|_{\cE}$ is a real analytic isomorphism on a neighbourhood of $0$ in $\cE$. {\it The latter means that an even function is uniquely determined by its spectrum.} There are still isospectral potentials to an even one, but they are not even.

On the other hand, given an arbitrary $q\in \cL^2[0,1]$, it is possible to determine it from its spectrum if additional data is considered. Let $\ell_1^2$ be defined as the set of sequences $(a_n)_{n\geq1} \in \ell^2(\R)$ such that
$$
\sum_{n\geq 1}n^2 a_n^2<\infty.
$$
Now, let $\kappa:\cL^2[0,1]\to \ell_1^2$ be defined by $\kappa(q)=(\kappa_n)_{n\geq1}$, with $\kappa_n=\log|y_2'(1,\lambda_n)|$. Recall that the eigenfunctions $y_2$ were defined at the beginning of this section as solutions of a specific initial value problem, the values $y_2'(1,\lambda_n)$ are called the \lq terminal velocities\rq.\  Each $\kappa_n$, is a real analytic function of $q$.

The sequence $\kappa(q)=(\kappa_n)_{n\geq1}$ changes sign when $q(x)\mapsto q(1-x)$ (this is a reflection over $x=1/2$), in contrast to $\Lambda(q)$ which is invariant under this reflection. Let $\Theta$ be a map defined by
\[
    \Theta: \cL^2[0,1] \to \ell_1^2\times \cS, \qquad
    \Theta(q)=(\kappa(q),\Lambda(q)).
\]
The fact that this map is injective on $\cL^2[0,1]$ and is a local real analytic isomorphism at every point in $\cL^2[0,1]$ is explained in the proof of Theorem 6 of Chapter 3 in \cite{PoeschTrub}.

Let us summarise the main results about isospectral potentials in the following theorem:
\begin{theorem}[see \cite{PoeschTrub}]
Let $p\in {\cL^2[0,1]}$, put
$$\cM(p)=\{q\in {\cL^2[0,1]}, \Spec(\cT_q) = \Spec(\cT_p) \}.$$
Then, $\cM(p)$ is a real analytic submanifold of {$\cL^2[0,1]$} that lies in the hyperplane consisting of all functions $q$ which mean satisfy $C(q)=C(p)$. Moreover, the set $\cM(p)$ is unbounded in {$\cL^2[0,1]$}.
\end{theorem}

The statement of this theorem consists of part $(a)$ of Theorem 1 and Theorem 2 of Chapter 4 in \cite{PoeschTrub}. Theorem 1 goes indeed further and gives a complete description of the tangent space to $\cM(p)$ and of its corresponding orthogonal complement. To conclude this section, we include the following observation. It connects the behaviour of the spectrum under the symmetries mentioned above with the corresponding well-known fact for Riemannian manifolds.

\begin{remark}
On a Riemannian manifold $M$, if $q\in \cC^{\infty}(M)$ and $\Psi$ is an isometry of $M$, then $\Delta + q$ and $\Delta + (q\circ \Psi)$ are isospectral. For the closed interval $[0,1]$, we have that $x\to 1-x$ is an isometry, therefore the potentials $q(x)$ and $p(x):=q(1-x)$ are isospectral, that is also proved in \cite{PoeschTrub} and this is why the spectrum does not distinguish left from right in this case.
\end{remark}

\section{Quasi-isospectral potentials, after Bilotta, Morassi, and Turco}
\label{sec:Quasi-Iso}

In this section, we present the construction of quasi-isospectral potentials, as is done in \cite{BILOTTA20171050}. Let $\cT_q$ be as it was defined in the previous section, and let  $\Spec(\cT_q) = (\lambda_m)_{m\geq 1}$ denote the spectrum of $\cT_q$. For each $m$, let $z_m$ denote the $\lambda_m$-eigenfunction normalised such that $z_m'(0) = 1$. In addition, $q$ should satisfy $\lambda_1(\cT_q) \geq 0$. Recall here that all the eigenvalues are simple.

The first question that arises here is whether there actually exist quasi-isospectral potentials. The answer is positive. In \cite{BILOTTA20171050},
Bilotta et al in  carry out an explicit construction of quasi-isospectral potentials, which we describe below.
The idea goes back to the problem of inserting or removing an eigenvalue of a given operator in a spectral gap.~

In Gesztesy-Telsch \cite{gesztesy-teschl}, the authors explain the single and the double commutation methods for general Sturm-Liouville operators on an interval, a line, and a half line.
The single commutation method allows us to introduce an eigenvalue inside a spectral gap, but has the disadvantage that it either introduces singularities or restricts the method to the use of strictly positive eigenfunctions. It goes back to the work of Jacobi \cite{jacobi}, Darboux \cite{darboux1915}, Crum \cite{crum}, Schmincke \cite{schmincke}, and Deift \cite{deift}. The double commutation method was introduced by Gel'fand and Levitan in 1955 \cite{gelfand_1951}, and it was improved by Gesztesy \cite{gesztesy} and Gesztesy and Teschl \cite{gesztesy-teschl}. More recently, Guliyev \cite{guliyev2020essentially}  investigated the properties of this method applied to isospectral Sturm-Liouville problems with boundary conditions involving rational Herglotz-Nevanlinna functions. The name of the method is related to the fact that it relies on a clever factorization of the operator and the commutation relation of its factors. Below, we will explain the general idea of this method. The approach used by Bilotta et al. has the same background idea. It consists of applying Darboux's lemma twice. Darboux's lemma, stated in Lemma~\ref{lemma:Darboux}, provides a new potential and a new eigenfunction. However, they both have singularities which are not square integrable on $[0,1]$. Despite that, in the same way as in the second commutation method, Bilotta et al. showed that using twice the lemma, one obtains an expression without singularities, in which all quantities are finite and well-defined. Thus, the method leads to regular solutions even though singular functions appear in the process. Let us explain here the parts of the construction that we consider relevant for this paper. The method is also described in \cite{PoeschTrub}.

\subsection{Construction of the potential and the eigenfunctions}\label{construction}

 Given a potential $q$, that we may consider in $\cL^2[0,1]$, let $(\lambda_m, z_m) _{m\geq 1}$ be a set of eigenvalues and eigenfunctions of $\cT_q$ satisfying $z'_m(0)=1$. Set $\lambda_0=0$, then they consider a fixed integer $n \geq 1$, and $t\in \R$ such that
\[ \lambda_{n-1} (q) < \lambda_{n} (q) + t < \lambda_{n+1} (q). \]
The construction in \cite{BILOTTA20171050} builds a potential $p$,
for which the spectrum of $\cT_p$ consists of the values
$$\lambda_m(q) + t\delta_{nm}, \quad \quad m\geq 1,$$
and whose eigenfunctions, denoted by $\kappa_{m,t}$, are computed explicitly in terms of various solutions of different initial value problems and boundary valued problems associated with the equation
    $$ -u''+qu=\lambda u.$$

Before we continue, let us state Darboux's Lemma, \cite[article 408]{darboux1915}:

\begin{lemma}[Darboux] \label{lemma:Darboux}
Let $\mu\in\mathbb{R}$ be fixed,  let $g$ be a non-trivial solution to
$$- u''+qu=\mu u.$$
Let $f$ be a non-trivial solution to
$$- u''+qu=\lambda u.$$
If $\lambda\neq\mu$,
denoting the Wronskian of $g$ and $f$ by $W[g,f]$, we have that
$$h=\dfrac{1}{g}W[g,f]:=\dfrac{1}{g}(gf'-g'f)$$
is a non-trivial solution to
\begin{equation}
	-u''+\left(q-2\dfrac{d^2}{dx^2}[\log(g)]\right)u=\lambda u.\label{Darboux_eq}
\end{equation}
If $\lambda=\mu$, the general solution to the last equation is then given by
$$\frac{1}{g} \left( a + b\int_0^x g^2(s) ds\right),$$
with $a, b\in \R$.
\end{lemma}
Notice that the function $h$, as well as the new potential, is defined on the complement of the zero set of the function $g$ in the interval $]0,1[$. Moreover, despite $g$ having simple zeros, $h$ would be only square integrable on compact sets in the complement of the zero set of $g$ in  $]0,1[$.

Following \cite{MORASSI2013288} and \cite{BILOTTA20171050}, consider  $y_1(x,q,\lambda_n +t), \ y_2(x,q,\lambda_n +t)$ such that for $x\in ]0,1[$
\begin{equation}
\left\{\begin{array}{l}
     -y_1'' + q y_1 = (\lambda_n + t) y_1 \\
     y_1(0) = 1,\quad y_1'(0)=0,
\end{array}\right.
\quad \quad
\left\{\begin{array}{l}
      -y_2'' + q y_2 = (\lambda_n + t)y_2 \\
           y_2(0) = 0,\quad y_2'(0)=1.
\end{array}\right. \label{problem_y1y2}
\end{equation}
Now set
$$w_{n,t}(x) := w(x,q,\lambda_{n}+t ):= y_1(x,q,\lambda_n +t) + \phi(t) y_2(x,q,\lambda_n +t)$$
with
$$\phi(t) = \frac{y_1(1,q,\lambda_n) - y_1(1,q,\lambda_n +t)}{y_2(1,q,\lambda_n +t)} .$$
Therefore, $w_{n,t}(x)$ satisfies
$$\left\{\begin{array}{l}
     -w_{n,t}'' + q w_{n,t} = (\lambda_n+t) w_{n,t} \\
     w_{n,t}(0)=1, \quad w_{n,t}(1)=y_1(1,q,\lambda_n).
\end{array}\right.$$

The first step is to apply Darboux's lemma to $w_{n,t}$, in the place of $f$, and in the place of $g$ use $z_n$.
Denoting
\begin{equation}
\beta_{n,t} = \beta_{n,t}(x,q) := W[w_{n,t},z_n] = w_{n,t} z'_n - w'_{n,t} z_n,
\label{eq:beta}
\end{equation}
and setting
$$h = \frac{\beta_{n,t}}{z_n},$$
We have that the latter satisfies
$$-h'' + \hat{q} h = (\lambda_n+t) h $$
with
$$\hat{q} = q - 2\frac{d^2\ }{dx^2} \left(\log{(z_n )}\right).$$
Notice that the function $h$ defined here is the opposite of the sum of the one given in Darboux's lemma. Since we are working with eigenfunctions, this change of sign is irrelevant. Additionally, notice here that both $h$ and $\hat{q}$ have singularities.
%%%%%%%%%%%%%%%%%%%%%%%%%%%%%%%%%%%%%%%%%%%%%%%%%%%%%%%%

The second step consists in applying the lemma to $z_m$ in the place of $f$ and $z_n$ in the place of $g$. In that way, we obtain a function
$\varphi$ that satisfies
$$-\varphi'' + \hat{q} \varphi = \lambda_m \varphi .$$

As the third step, we apply the lemma to $\varphi$ in the place of $f$ and $h$ in the place of $g$. In this way, we obtain a function
$\kappa_{m,t}$
that satisfies
$$-\kappa_{m,t}'' + p \ \kappa_{m,t} =  \lambda_m \kappa_{m,t}$$
with
\begin{equation}
p= \hat{q} - 2\frac{d^2\ }{dx^2} \left(\log{(h)}\right)= q - 2\frac{d^2\ }{dx^2} \left(\log{(z_n h)}\right) = q - 2\frac{d^2\ }{dx^2} \left(\log{(\beta_{n,t})}\right).
\label{eq:p}
\end{equation}

A key observation here is that if $q\in \cL^2[0,1]$, the function $\beta_{n,t}(x,q)$ is continuous and strictly positive for
$$(x,t)\in [0,1]\times ]\lambda_{n-1}-\lambda_n, \lambda_{n+1}-\lambda_n[.$$
This fact requires a proof that is given in \cite[Lemma 1 Ch.6]{PoeschTrub}.

We now consider the cases $m\neq n$ and $m=n$ separately.

\textit{Case 1: $m\neq n$}. We consider $m\in \N$, $m\geq 1$ and $m\neq n$. In this case, the functions $\varphi$ and $\kappa_{m,t}$ defined above and determined by Darboux's lemma are given by
$$\varphi = \frac{1}{z_n} W[z_m, z_n] \quad \text{ and } \quad \kappa_{m,t} = \frac{1}{h} W[\varphi, h],$$
respectively. Another key observation at this point is that
$$\kappa_{m,t}(x)= (\lambda_n-\lambda_m)z_m(x)-\dfrac{1}{z_n(x)}[z_n(x),z_m(x)]\left(\dfrac{d}{dx}[\log(\beta_{n,t}(x))]\right),$$
see \cite[Lemma 2 Ch.5]{PoeschTrub}. After this statement, it is easy to see that
$$\kappa_{m,t}(x) = (\lambda_n-\lambda_m)\left(z_m(x) - \frac{t \, w_{n,t}}{\beta_{n,t}} \int_{0}^{x} z_n(s) z_m(s) \, ds \right),$$
from which is clear that $\kappa_{m,t}(x)$ is well-defined and has no singularities in $\cL^2[0,1]$.

\textit{Case 2: $m= n$}. Setting
$$\kappa_{n,t}(x) = \dfrac{z_n(x)}{\beta_{n,t}(x)}.$$
A direct computation shows that $\kappa_{n,t}$ satisfies
$$-\kappa_{n,t}'' + p \ \kappa_{n,t} =  (\lambda_n+t) \kappa_{n,t}.$$
with $p$ as in equation (\ref{eq:p}), as well as the Dirichlet boundary conditions,
$$\kappa_{n,t}(0) = \kappa_{n,t}(1) = 0.$$

Moreover, $\kappa_{m,t}$ also satisfies the Dirichlet condition since so does $z_m$ and both $z_m'(0)$ and $z_{m}'(1)$ are non-zero and this implies that
$$\frac{1}{z_n(x)}[z_n(x),z_m(x)]\to 0 \quad \text{as} \quad  x\to 0 \quad  \text{or} \quad  x\to 1.$$
Let us summarize the previous results (from \cite{BILOTTA20171050, PoeschTrub}) into the following lemma.
\begin{lemma}[Bilotta et al.] Using the notation defined above, the functions
$$\kappa_{m,t}(x)=\left\{\begin{array}{cc}
    (\lambda_n-\lambda_m)z_m(x)-\dfrac{1}{z_n(x)}[z_n(x),z_m(x)]\left(\dfrac{d}{dx}[\log(\beta_{n,t}(x))]\right) & \text{ if } m\neq n  \\
    \dfrac{z_n(x)}{\beta_{n,t}(x)} & \text{ if } m=n
\end{array}\right.$$
are well-defined and are non-trivial solution of the problem
$$u''+\left({q}-2\dfrac{d^2}{dx^2}[\log(\beta_{n,t})]\right)u=(\lambda_m+t\delta_{nm}) u,$$
satisfying the Dirichlet boundary conditions. Furthermore, if $q\in \cL^2[0,1]$, $\beta_{n,t}$ is strictly positive, $\beta_{n,t}^{''} \in \cL^{2}[0,1]$, and $p \in \cL^{2}[0,1]$.
\end{lemma}
The last statements about $\beta_{n,t}$  and $p$ follow from the definition of $\beta_{n,t}$ in equation (\ref{eq:beta}) which implies that
$$\beta_{n,t}^{''} = t (z_{n}'w_{n,t}-z_{n}w_{n,t}').$$
\begin{corollary} \label{c:qipSLo}
With the notation introduced above, we have that the potentials $q$ and $p = {q}-2\tfrac{d^2}{dx^2}[\log(\beta_{n,t})]$ are quasi-isospectral.
\end{corollary}
\subsection{The boundary conditions in Darboux's lemma}

The result in this section follows from a straightforward computation; we include it here for completeness, as it may not be explicitly stated in the literature.
Let us consider Darboux's lemma as stated in Lemma~\ref{lemma:Darboux}. First of all, it is noteworthy that Darboux's lemma does not mention any boundary conditions.

\begin{lemma}
Let $f$ and $g$ satisfy the following mixed problems,
$$\left\{\begin{array}{l}
-f''+qf=\lambda_f f\\
\alpha_1 f(0)+ \alpha_2 f'(0)=0\\
\beta_1 f(1)+\beta_2 f'(1)=0
\end{array}\right.  \quad \text{ and } \quad
\left\{\begin{array}{l}
-g''+qg=\lambda_g g \\
\tilde{\alpha}_1 g(0) + \tilde{\alpha}_2 g'(0)=0\\
\tilde{\beta}_1 g(1) + \tilde{\beta}_2 g'(1)=0
\end{array}\right. $$
respectively, such that $\alpha_{1} \cdot \alpha_{2} \cdot \beta_{1} \cdot \beta_{2}\neq 0$, and $\tilde{\alpha}_{1} \cdot \tilde{\alpha}_{2} \cdot \tilde{\beta}_{1} \cdot \tilde{\beta}_{2}\neq 0$.

Then, the non-trivial solution $h:= \frac{1}{g} W\left[g,f\right]$ of $$-u''+\tilde{ q } \ u=\lambda_f \ u,$$
given by Darboux's lemma, satisfy the Dirichlet boundary conditions if
$$\tilde{ \alpha_{1}}\alpha_{2}-\alpha_{1}\tilde{ \alpha_{2} }=\tilde{ \beta_{1}}\beta_{2}-\beta_{1}\tilde{ \beta_{2} }=0.$$
In particular, this applies if $f$ and  $g$ satisfy the same non-Dirichlet boundary conditions.

On the other hand, if $f$ and $g$ are solutions of the Dirichlet problem, $h$ is not defined at $x=0$, but the following limits vanish:
$$\lim_{x\to 0^+} h(x)  = 0 = \lim_{ x\to 1^- }h(x).$$
\end{lemma}

\begin{proof}
Using the definition of $h$ and the conditions that $f$ and $g$ satisfy, we obtain that the boundary conditions for $h$ are given by
$$\lim_{x\to 0} h(x)=\lim_{ x\to 0 } \left(\dfrac{\tilde{\alpha_1} }{\tilde{ \alpha_2 } }-\dfrac{{\alpha_1} }{{ \alpha_2 } } \right)f(x),\quad  \lim_{ x\to 1 }h(x)=\lim_{ x\to 1 } \left(\dfrac{\tilde{\beta_1} }{\tilde{ \beta_2 } }-\dfrac{{\beta_1} }{{ \beta_2 } } \right)f(x).$$

Dirichlet boundary conditions for $h$ will then be obtained if
$$\tilde{ \alpha_{1}}\alpha_{2}-\alpha_{1}\tilde{ \alpha_{2} }=\tilde{ \beta_{1}}\beta_{2}-\beta_{1}\tilde{ \beta_{2} }=0.$$
This applies, in particular,  if $f$ and $g$ satisfy the same boundary conditions that are not Dirichlet.

On the other hand, if we consider the Dirichlet problem, $h$ is not defined at $x=0$, but we can compute the limit, which exists because the zeros of $g$ and of $f$ are simple. So, in the same way as in the previous section, we have that
\begin{eqnarray*}
\lim_{x\to 0^+} h(x) &=& f'(0) - \lim_{ x\to 0^+ } \left(\dfrac{f(x) }{g(x)}\right)g'(0)
= f'(0) - \left(\dfrac{f'(0)}{g'(0)} \right)g'(0) = 0,\\
\lim_{ x\to 1^- }h(x)&=& f'(1) - \lim_{ x\to 1^- } \left(\dfrac{f(x) }{g(x)}\right)g'(1)
= f'(1) - \left(\dfrac{f'(1)}{g'(1)} \right)g'(1) = 0.
\end{eqnarray*}

\end{proof}

\subsection{The effect of Darboux's lemma on the spectrum}

The result presented in this section appears to be of some independent interest; to our knowledge, it has not been explicitly stated in the literature.
It consists of some observations about the effect on the spectrum of Schr\"odinger-type operators on $ [0,1] $ of applying the commutation method and Darboux's lemma.

Applying Lemma~\ref{lemma:Darboux} can lead to the removal of an eigenvalue. \[
T_g := g\left(\frac{d}{dx}\right)\frac{1}{g}
\]
on {$ \mathcal{L}^{2}[0,1] $}. Its formal adjoint is given by
\[
T_{g}^* = - \frac{1}{g} \left(\frac{d}{dx}\right) g,
\]
and these operators yield to the factorizations
\[
\mathcal{T}_q - \mu  = T_g^* T_g, \qquad \mathcal{T}_{\tilde{q}} - \mu  = T_g T_g^*,
\]
where $\tilde{q}$ denotes the new potential in (\ref{Darboux_eq}).
From the commutation formula (Deift \cite[Theorem~2]{deift}), it follows that
\[
\operatorname{Spec}(\mathcal{T}_q - \mu ) \setminus \{0\} = \operatorname{Spec}(\mathcal{T}_{\tilde{q}} - \mu) \setminus \{0\}.
\]
However, it remains unclear whether $0$ is an eigenvalue of $\mathcal{T}_{\tilde{q}} - \mu $. To explore this, assume $0 \in \operatorname{Spec}(\mathcal{T}_{\tilde{q}} - \mu )$. Then, there exists $\phi \in {\mathcal{L}^{2}[0,1]}$ such that $ (\mathcal{T}_{\tilde{q}} - \mu )(\phi) = 0 $, and it follows that
\[
\left\| T_g^* \phi \right\|^2 =
\left\|\frac{1}{g} \frac{d}{dx}(g \phi)\right\|^2 = \langle \phi, (\mathcal{T}_{\tilde{q}} - \lambda)(\phi) \rangle = 0.
\]
Therefore, $g\phi$ should be constant, and $\phi$ should be a scalar multiple of {$\frac{1}{g}$}.

For example, if we consider the Dirichlet problem, and if $g = z_m$ is the $m$-th Dirichlet eigenfunction, then $g$ has exactly $m+1$ simple zeros in $[0,1]$. Then, $\frac{1}{g} \notin {\mathcal{L}^2[0,1]}$, which
implies that $0 \notin \operatorname{Spec}(\mathcal{T}_{\tilde{q}} - \mu)$.
This means that, in this example, the application of Darboux's lemma removes the eigenvalue $ \mu $ from the spectrum.

Let us now look at the second part of Lemma~\ref{lemma:Darboux}, which addresses the case when the eigenvalues corresponding to $f$ and $g$ are the same, i.e. $\mu = \lambda$. The general solution to the equation $(\mathcal{T}_{\tilde{q}} - \lambda)\phi = 0$ is given by
\[
h := \frac{1}{g}\left( a + b \int_0^x g^2(s)\, ds \right),
\]
where \( a, b \in \mathbb{R} \). Let us now define a new operator on ${\mathcal{L}^2[0,1]}$:
\[
T_{h} = h \left(\frac{d}{dx}\right) \frac{1}{h},
\]
whose formal adjoint is
\[
T_{h}^* = - \frac{1}{h} \left(\frac{d}{dx}\right) h.
\]
This gives rise to the factorization:
\[
T_{h}^* T_{h} = T_g T_g^* = \mathcal{T}_{\tilde{q}} - \lambda, \quad \text{and} \quad T_{h} T_{h}^* = \mathcal{T}_{p} - \lambda.
\]
By Deift’s theorem, we have
\[
\operatorname{Spec}(\mathcal{T}_{\tilde{q}} - \lambda) = \operatorname{Spec}(\mathcal{T}_{p} - \lambda) \setminus \{0\}.
\]
In contrast to the earlier case, $0$ is an eigenvalue of $\mathcal{T}_{p} - \lambda$. Indeed, if $\phi \in {\mathcal{L}^2[0,1]}$ satisfies $(\mathcal{T}_{p} - \lambda)\phi = 0$, then $\phi$ is a scalar multiple of $\frac{1}{h}$. Unlike $\frac{1}{g}$, the function $\frac{1}{h}$ belongs to {$\mathcal{L}^2[0,1]$} for a suitable choice of constants $a$ and $b$. Hence,
\[
\operatorname{Spec}(\mathcal{T}_q - \lambda) = \operatorname{Spec}(\mathcal{T}_{p} - \lambda).
\]
Moreover, we obtain the following explicit expression of the new problem:
\begin{align*}
T_{h} T_{h}^*&= -\frac{d^2}{dx^2} + \left[q - 2 \frac{d^2}{dx^2} \log\left(a + b \int_0^x g^2(s)\, ds\right) \right] - \lambda.
\end{align*}
Consequently, the singularities are removed, and the corresponding solutions are well-behaved on the interval $[0,1]$.

\section{Examples}
\label{sec:Examples}
In this section, we present several examples illustrating different aspects of the behavior of quasi-isospectral potentials. To the best of our knowledge, these examples have not been previously discussed in the literature.

\subsection{Example of a family of even quasi-isospectral potentials}

In this section, we construct a family of even potentials that are quasi-isospectral to $q = 0$. This family is parametrised by the positive eigenvalues of the Dirichlet Laplacian on the unit interval and by the choice of the new eigenvalue ($\mu_n$ below).
Let us start putting $q = 0$, so $\mathcal{T}_q = \mathcal{T}_0 = \Delta_{\mathcal{I_D}}$. The spectrum of this operator is well known, and it was described in Example~\ref{ex:eigenspace}. We are going to use the notation introduced in Section~\ref{construction}. In this construction, only one eigenvalue is perturbed.
More precisely, if
\[
    \Spec(\mathcal{T}_0)
    = (\lambda_k)_{k \geq 1}
    = (k^2\pi^2)_{k \geq 1}.
\]
Let $\mu_n \in \mathbb{R}$ be such that
\[
(n-1)^2 \pi^2< \mu_n < (n + 1)^2 \pi^2, \quad (n \geq 2).
\]
In Section~\ref{sec:StLiD}, we pointed out that an even function is uniquely determined by its spectrum.
Our goal in this section is to use Corollary~\ref{c:qipSLo} to explicitly describe the unique even potential $p$ that satisfies $ \text{Spec}(\mathcal{T}_{p})= \left(\widetilde{\lambda}_{k}\right)_{k \geq 1}$, where
$$ \widetilde{\lambda}_{k}=\lambda _{k}  \quad \text{ for all } \quad k \neq n, \quad \text{ and }  \quad \widetilde{\lambda}_{n}=\mu _{n} .$$
The solutions to the initial value Problems~\ref{problem_y1y2} are well known and given explicitly by:
\[
y_1(x, \lambda) = \cos(\sqrt{\lambda} x), \quad y_2(x, \lambda) = \frac{\sin(\sqrt{\lambda} x)}{\sqrt{\lambda}},
\]
where $\lambda$ here corresponds to $\lambda_n + t = \mu_n$. By construction, $z_n(x)= y_2(x, \lambda_n)$ and $w_n$ is determined by the linear combination:
\begin{align*}
    w_n(x, \lambda) &= \cos(\sqrt{\lambda} x) + \frac{\cos(\sqrt{\lambda_n}) - \cos(\sqrt{\lambda})}{\sin(\sqrt{\lambda})} \sin(\sqrt{\lambda} x).
    \intertext{Using trigonometric identities:
	    \begin{align*}
		    \cos(a x)\sin(a) &= \tfrac{1}{2} \left( \sin(a(1 - x)) + \sin(a(1 + x)) \right), \\
		    \cos(a)\sin(a x) &= \tfrac{1}{2} \left( \sin(a(1 + x)) - \sin(a(1 - x)) \right),
	    \end{align*}
    we can simplify $w_n$ as follows:}
               w_n(x, \lambda) &= \frac{\cos(\sqrt{\lambda} x) \sin(\sqrt{\lambda}) + (\cos(n\pi) - \cos(\sqrt{\lambda})) \sin(\sqrt{\lambda} x)}{\sin(\sqrt{\lambda})}\\
    &= \frac{\sin(\sqrt{\lambda}(1 - x)) + (-1)^n \sin(\sqrt{\lambda} x)}{\sin(\sqrt{\lambda})}.
\end{align*}

Differentiating, we obtain:
\begin{align*}
    w_n'(x, \lambda) &= \frac{-\sqrt{\lambda} \cos(\sqrt{\lambda}(1 - x)) + (-1)^n \sqrt{\lambda} \cos(\sqrt{\lambda} x)}{\sin(\sqrt{\lambda})}.
\end{align*}

From the expression for the function $\beta_{n,t}(x)$ in (\ref{eq:beta}), we have:
\begin{align*}
        \beta_{n,t}(x) &= \frac{ (\sin(\sqrt{\mu_n}(1 - x)) + (-1)^n \sin(\sqrt{\mu_n} x)) \cos(n\pi x) }{ \sin(\sqrt{\mu_n}) } \\
    &\quad + \frac{ (\sqrt{\mu_n} \cos(\sqrt{\mu_n}(1 - x)) + (-1)^{n+1} \sqrt{\mu_n} \cos(\sqrt{\mu_n} x)) \sin(n\pi x) }{ n\pi \sin(\sqrt{\mu_n}) }.
\end{align*}

Using the symmetry properties of the trigonometric functions:
\[
\cos(n\pi(1 - x)) = (-1)^n \cos(n\pi x), \quad \sin(n\pi(1 - x)) = (-1)^{n+1} \sin(n\pi x),
\]
we deduce that
\[
\beta_{n,t}(1 - x) = \beta_{n,t}(x), \quad \text{for all } x \in [0,1].
\]

This shows that $\beta_{n,t}$ is an even function, and consequently, the potential
$$ p = -2 \frac{d^2\ }{dx^2}\log\left(\beta_{n,t}\right) = -2 \left( \frac{\beta_{n,t}^{''}}{\beta_{n,t}} - \left(\frac{\beta'_{n,t}}{\beta_{n,t}}\right)^2 \right) $$
is also even. In this way, we have constructed a whole family of even potentials that, by Corollary~\ref{c:qipSLo} are quasi-isospectral to $q = 0$.

\subsection{Examples with different boundary conditions}

A natural question that arises is whether two operators defined by the same Sturm-Liouville differential expression but with different boundary conditions may be isospectral or quasi-isospectral. This question traces back to the foundational work of Sturm \cite{sturm1836equations} and Liouville \cite{liouville1836memoire} in $1836$, who first studied how boundary conditions influence the spectrum of such operators. We start with a simple example to throw light on this phenomenon. These are standard computations that we include here only for the sake of illustration.

\begin{example}
Let us consider the following Sturm-Liouville problem
\begin{equation}
\begin{cases}
- y'' (x) + q(x) y(x) = \lambda y(x), & x\in ]0,1[ \\
a(t) y(0) + b(t) y'(0) = 0, & \\
c(t) y(1) + d(t) y'(1) = 0. &
\end{cases} \label{eq:SLp-vbc}
\end{equation}
where $a(t), b(t), c(t), d(t)$ are functions depending on a parameter $t$. We want to find conditions on these functions that yield a one-parametric family of isospectral operators, which are closed extensions of the initial differential operator whose minimal domain is $\cC_c^{\infty}[0,1]$.

{\bf Claim:} If $q(x) =0$, boundary conditions satisfying
\[ a(t)d(t) - b(t) c(t) =0 \]
yield a solution to the problem in (\ref{eq:SLp-vbc}) when $\lambda = \pi^2 n^2 \neq 0$.

Let $q(x) =0$. The general solution for $y''(x) -\lambda y(x)=0$ has the form
$$y(x)=c_1\cos{(\sqrt{\lambda}x)}+c_2\sin{(\sqrt{\lambda}x)}.$$
Thus, $y(0)=c_1$, $y'(0)=c_2\sqrt{\lambda}$ and the first boundary condition gives:
\begin{equation*}
    c_1 a(t) +c_2 \sqrt{\lambda} b(t) = 0 .
    \end{equation*}
Moreover,
\begin{align*}
    y(1)&=c_1\cos{(\sqrt{\lambda})}+c_2\sin{(\sqrt{\lambda})},\\
    y'(1)&=-c_1\sqrt{\lambda}\sin{(\sqrt{\lambda})}+c_2\sqrt{\lambda}\cos{(\sqrt{\lambda})},
\end{align*}
and the second boundary condition becomes:
$$c_1\left(c(t) \cos{(\sqrt{\lambda})}- d(t)\sqrt{\lambda}\sin{(\sqrt{\lambda})}\right)+c_2 \left(c(t) \sin{(\sqrt{\lambda})}+ d(t) \sqrt{\lambda}\cos{(\sqrt{\lambda})}\right)=0.$$

In order to have a non-trivial solution, the system should satisfy:
\[ a(t)  \left(c(t) \sin{(\sqrt{\lambda})}+ d(t) \sqrt{\lambda}\cos{(\sqrt{\lambda})}\right) = \sqrt{\lambda} b(t) \left(c(t) \cos{(\sqrt{\lambda})}- d(t)\sqrt{\lambda}\sin{(\sqrt{\lambda})}\right) \]
that can be rewritten as
\begin{equation} \left(  a(t) c(t) + b(t) d(t) \lambda \right) \sin{(\sqrt{\lambda})} + \left( a(t)d(t) - b(t) c(t) \right) \sqrt{\lambda} \cos{(\sqrt{\lambda})} = 0.
\label{eq:ceevdbc}
\end{equation}

If $\sqrt{\lambda}=n\pi$ we obtain:
\[ a(t)d(t) - b(t) c(t) =0. \]
A condition that makes sense is $a(t) = c(t)$ and $b(t) = d(t)$. Therefore, the boundary conditions should have the form:
$$\begin{cases}
a(t) y(0) + b(t) y'(0) = 0, & \\
a(t) y(1) + b(t) y'(1) = 0. &
\end{cases}
$$
Conversely, setting $a(t)d(t) - b(t) c(t) =0$ in equation (\ref{eq:ceevdbc}), one obtains
$$ \left( a(t) c(t) + b(t) d(t) \lambda \right) \sin{(\sqrt{\lambda})} =0$$
which occurs when $\sqrt{\lambda}=n\pi$.

On the other hand, if we want to interpolate between Dirichlet and Neumann boundary conditions, it makes sense to choose
$$a(t) = \cos{(t)}, \  b(t)=\sin{(t)}, \  \ t\in [0,\pi/2].$$
This is the type of boundary conditions considered by Sturm and by Liouville. However, notice that $0\in \Spec(\Delta_N)$ but $0\notin \Spec(\Delta_D)$. In addition, $\lambda=0$ is not an eigenvalue of the corresponding extension if $t >0$. So, in this case, we have a one-parameter family of isospectral extensions if $t\in [0,\pi/2[$. If the potential $q$ is constant, the spectrum just gets shifted, and we have a similar result for constant potentials.
\end{example}

\begin{example}
In this example, we apply Darboux's lemma to construct two isospectral potentials whose eigenfunctions obey different boundary conditions.
For any $q\in {\cL^2[0,1]}$, let us consider the following Sturm Liouville problem
$$\begin{cases}
- y'' (x) + q(x) y(x) = \lambda y(x), & x\in ]0,1[ \\
y'(0) - h y(0) = 0, & h\in\mathbb{R} \\
y'(1) + H y(1) = 0, & H\in\mathbb{R}
\end{cases}
$$
As in previous cases, this problem has a discrete spectrum $(\lambda_k)_{k\geq1}$, each eigenvalue being simple, \cite{Freiling2001InverseSP}. Let $(\phi_k)_{k\geq1}$ be the corresponding eigenfunctions. Applying the lemma of Darboux to $g=\phi_1$ and $f=\phi_k$, for each $k\geq2$, we obtain $h_k=\frac{1}{\phi_1}W[\phi_1,\phi_k]$ is a non-trivial solution of
$$-u''+\tilde{q}u=\lambda_k u,$$ where
$\tilde{q}=\left(q-2\dfrac{d^2}{dx^2}[\log(\phi_1)]\right)$
Moreover, $h_k$ is a solution to the Dirichlet problem:
$$\begin{cases}
- y'' (x) + \tilde{q}(x) y(x) = \lambda_k y(x), & x\in ]0,1[ \\
y(0)=0=y(1).
\end{cases}
$$
We obtain the same result if we apply the lemma to $g=\phi_1$ and $f=\phi_1$, but with a different $h_1$. Thus, we have obtained that $q$ and $\tilde{q}$ are isospectral with different boundary conditions.
\end{example}

\section{Compactness of isospectral manifolds and potentials}
\label{s:cimap}
In this expository section, we summarize and briefly mention some known results about compactness of isospectral Riemannian manifolds and of isospectral potentials.

\subsection{Compactness of isospectral Riemannian manifolds}
Compactness and precompactness results for isospectral manifolds have been established in various settings. Let us briefly mention some key results and contributors to this area.

One of the earliest results is due to Richard Melrose, who showed in \cite{Mel83}  that isospectral sets of drum-shaped domains are precompact in the smooth topology. Later, Osgood, Phillips, and Sarnak proved in \cite{OPS1} that isospectral sets of smooth metrics on compact surfaces are compact in the $C^\infty$-topology (modulo isometries). Key tools in their proof are heat invariants and the regularised determinant of the Laplacian.

Numerous extensions of the compactness results by Osgood, Phillips, and Sarnak have been established in various geometric settings. For instance, Kim \cite{Kim08} extended their work to the case of flat surfaces with a boundary. In higher dimensions, compactness of isospectral metrics within a fixed conformal class has been studied by Brooks, Perry, and Yang \cite{Brooks-Perry-Yang}, Chang and Yang \cite{CY90}, and Chang and Qing \cite{CQ97}. Additional compactness results under geometric constraints such as uniform bounds on curvature or injectivity radius have been obtained by Brooks and Glezen \cite{bg94}, Brooks, Perry, and Petersen \cite{BPP92},
and Zhou \cite{Zhou}. For an overview of both positive and negative results in the isospectral setting, we refer the reader to the survey by Gordon, Perry, and Schueth \cite{GPS}.

In the noncompact setting, where traditional spectral invariants may be insufficient, one considers quantities like the scattering phase or resonances, and speaks of \emph{iso-phasal} or \emph{iso-resonant} families. Hassell and Zelditch \cite{HZ99} established compactness of iso-phasal exterior domains in the plane, in a natural smooth topology. Additional contributions are the results by Borthwick and Perry \cite{BP11}, and
Borthwick, Judge, and Perry \cite{BJP03} about compactness of open hyperbolic manifolds with funnel ends.

The first author of this paper proved in \cite{Ald10} the precompactness of iso-resonant hyperbolic surfaces with cusps under additional conditions.  Later, in joint work with Pierre Albin and Fr\'ed\'eric Rochon, showed compactness of relatively isospectral surfaces with funnels, cusps, and boundaries, \cite{AARrivcs}.

The notion of compactness is often obtained either via Sobolev embeddings and Re\-llich's compactness theorem, or through compactness theorems specifically designed for manifolds or metric spaces. One of the best-known results in this area is Gromov’s compactness theorem, which provides conditions -namely, uniform bounds on diameter and on the number of balls of a fixed radius covering the space under which a sequence of metric spaces admits a convergent subsequence in the Gromov–Hausdorff topology. Other compactness results depend on the particular nature of the problem and the geometric quantities one seeks to control. For example, the Cheeger–Gromov compactness theorem gives sufficient conditions to extract a smoothly convergent subsequence of Riemannian manifolds, assuming uniform bounds on geometric quantities such as Ricci curvature, scalar curvature, diameter, and injectivity radius.

\subsection{Compactness of isospectral potentials.}
\label{subs:cisp}
In this section, we recall the compactness results for isospectral potentials on closed manifolds by Br\"uning \cite{Bruning1984OnTC} and Donnelly \cite{Donnelly2004CompactnessOI}, as well as the more recent result of compactness of isospectral potentials for bounded sets of $\R^n$ under certain conditions by Choulli et al \cite{Choulli2015HeatTA}.

Let $(M,g)$ be a closed Riemannian manifold of dimension $n$ and $\Delta$ be the corresponding (positive) Laplace-Beltrami operator.  Let $V\in \cC^{\infty}(M)$, and denote $H=\Delta+V$ the Schr\"odinger operator associated to $V$. Let $K_H$ denote the heat kernel corresponding to the operator $H$. As we already mentioned, the heat kernel has an asymptotic expansion as $t \to 0^+$ whose restriction to the diagonal can be written as
$$K_H(t,x,x)=\dfrac{1}{(4\pi t)^{n/2}}(1+ u_1(x)t + u_2(x)t^2 + \cdots)$$
and the terms $u_i(x)$ are invariant polynomials in the components of the curvature tensor on $M$, the potential function $V$, and their covariant derivatives of higher order, see \cite{Gilkey1979RecursionRA}. If $\Spec(H)=(\mu_k)_{k\geq1}$, then as $t\to 0^+$

\begin{equation}
\text{Tr}\left(e^{-Ht}\right)=\sum_{k=1}^{\infty}e^{-\mu_k t}\sim \dfrac{1}{(4\pi t)^{n/2}}\left(\text{Vol}(M)+ a_1 t + a_2 t^2+\cdots\right) \label{eq:aeHtSoSt}\end{equation}
where $a_1=\int_M [V+ v_{11}] d\mu $, $a_2=\frac{1}{2}\int_M [V^2+ v_{21}V + v_{22}]d\mu$, and for all $j\geq3$
$$a_j=\int_M u_j d\mu = c_j\int_M |\nabla^{j-2}V|^2+\sum_{k=0}^j\left(\sum_{\alpha_k \in \Z^k, |\alpha_k| \leq\ell(k)}\int_M P_{\alpha_1}^k(V)P_{\alpha_2}^k(V)\cdots P_{\alpha_k}^k(V)\right) d\mu$$
where the functions $v_{ij}$ are smooth on $M$ depending only on the metric $g$ but not on the potential $V$.  Each $P_{\alpha_i}^k$ is a differential operator with smooth coefficients depending only on the metric. Furthermore, $\text{ord}P^k_{\alpha_i} \leq j-3$ and $\sum_{i=1}^{k}\text{ord}P^k_{\alpha_i} \leq 2(j-3)$ see \cite[Prop.3.1 and Prop.3.2]{Donnelly2004CompactnessOI} and \cite{Bruning1984OnTC}.

We define the set of potentials isospectral to $V$ as
$$\text{Is}(V) = \{W \in \cC^{\infty}(M) \ | \ \Spec(\Delta + V) = \Spec(\Delta + W),\  \text{ including multiplicities}\}.$$
Notice that the heat invariants of isospectral potentials are the same. Now,
define
$$n_0=\inf\{m\in \N | m\geq n/2 \}.$$

Let $W_{s,r}(M)$ denote the usual $(s,r)$ Sobolev space on $M$. If $s\in \N$, $W_{s,r}(M)$ consists of functions in $\cL^{r}(M)$ that have $s$ weak derivatives in $\cL^{r}(M)$.
Then, using the Sobolev embedding theorems, Br\"uning and Donnelly respectively obtained:
\begin{theorem}[Br\"uning \cite{Bruning1984OnTC} Thm.2]
If $\text{Is}(V)$ is bounded in $W_{3n_0-2,2}(M)$ then it is compact in $\cC^\infty(M)$. This implies that if $\dim(M)\leq3$, then $\text{Is}(V)$ is compact in $\cC^\infty(M)$.
\end{theorem}

\begin{theorem}[Donnelly \cite{Donnelly2004CompactnessOI} Thm.4.1] If $\text{Is}(V)$ is bounded in $W_{s,2}(M)$ with $s > n/2-2$ then it is compact in $\cC^\infty(M)$. This implies that if $\dim(M)\leq 9$, $\text{Is}^+(V)$, the set of non-negative potentials isospectral to $V$, is compact in $\cC^\infty(M)$. \label{DonnellyCompact}
\end{theorem}

The idea of the proof is to show that the heat invariants corresponding to isospectral potentials are equal. Then, one uses the description of the heat invariants given above to obtain uniform estimates in Sobolev norms of the form
$$\|W\|_{j,2} \leq C_{j,g}, \quad \text{ for all}\quad W\in \text{Is}(V).$$
Compactness then follows from the Sobolev embeddings. Additionally, as Br\"uning pointed out, thanks to the Gagliardo-Nirenberg inequalities, the uniform estimate given in the statement of the theorem, instead of the uniform estimates in all the Sobolev norms, is sufficient to obtain the result.

Now, let $\Omega\subset\mathbb{R}^n$, for $n\geq1$, be a bounded connected domain with smooth boundary. Let $A$ be a second-order elliptic differential operator of Laplacian type with domain $H^1_0(\Omega)$. The Dirichlet Laplacian on $\Omega$ is a particular case of $A$. Let $V\in \cC_0^\infty (\Omega)$ be a real-valued potential, and $H:=A+V$. Then, the results of Choulli et al. in \cite{Choulli2015HeatTA} imply that, taking $A$ as the Dirichlet Laplacian $\Delta$, there exists an asymptotic expansion of the relative heat trace as $t\to 0^+$

\begin{equation}
\Tr \left(e^{-tH}-e^{-t \Delta} \right) \sim t^{-n/2} \left(a_1(V) t + \cdots+a_i(V) t^i + O(t^{i+1})\right)
\label{eq:ae-rht-st-sch}
\end{equation}
where the $a_j(V)$ are coefficients. In addition, since each of the heat operators is trace class, one has that
$$\Tr \left(e^{-tH}\right)=\Tr \left(e^{-t\Delta} \right)+\Tr \left(e^{-tH}-e^{-t\Delta} \right).$$
Therefore, from calculations involving the first four terms in the asymptotic expansion, they obtain:
\begin{theorem}[Choulli et al. \cite{Choulli2015HeatTA} Thm. 1.1]
For any bounded subset $\mathcal{B}\subset L^{\infty}(\Omega)$ such that $V\in\mathcal{B}$, the set $\mathcal{B}\cap \text{Is}(V)$ is compact in $H^s(\Omega)$ for each $s\in]-\infty,2[$.
\end{theorem}

\section{Quasi-isospectral Riemannian manifolds}
\label{s:qiRm}
This section contains the main new contributions of this paper.
Let $(M,g)$ be a closed connected Riemannian manifold. Let $\Delta_g$ be the associated Laplace operator. Let us first define quasi-isospectral Riemannian manifolds.

\begin{definition} Two Riemannian manifolds $(M_1,g_1)$, $(M_2,g_2)$ are called quasi-isospectral  if the spectra of the corresponding Laplace operators are the same with the exception of one eigenvalue that is allowed to change in a prescribed interval. More precisely, if $k\geq 1$ and $\varepsilon >0$, we say that they are $(k,\varepsilon)$-quasi-isospectral if $\lambda_{m}(g_1) = \lambda_{m}(g_2)$ for $m\neq k$, and
$$\lambda_k(g_1) \in ]\lambda_k(g_2)- \varepsilon, \lambda_k(g_2) + \varepsilon[.$$
Moreover, we require that
$$\lambda_{k-1}(g_2) < \lambda_k(g_2) - \varepsilon \quad \text{ and }  \quad \lambda_k(g_2) + \varepsilon < \lambda_{k+1}(g_2).$$
We can drop $k$ in the notation to say that the metrics are $\varepsilon$-quasi-isospectral if we do not need to indicate which eigenvalue changes.
\end{definition}

Notice here that $\lambda_0(g)=0\in \Spec(\Delta_g)$, here we only consider perturbations of positive eigenvalues.

On a closed manifold, the heat operator is trace class, and the heat trace has an asymptotic expansion for small values of $t$  of the form
\begin{equation}
\Tr\left(e^{-t\Delta_g}\right) = t^{-n/2} \sum_{j=0}^{\infty} a_j(g) t^j \label{eq:aehtLnb}
\end{equation}
where $n$ is the dimension of $M$ and the coefficients $a_j(g)$ are the heat invariants that depend on the metric and its derivatives in a precise form. We know that isospectral metrics have the same heat invariants. What can we say about quasi-isospectral metrics? The following theorem gives as an answer.

\begin{theorem}
Let $M$ be a closed manifold with $n=\dim(M)$. Let $g$ be a fixed Riemannian metric on $M$, and define
$$QI_\varepsilon(g) = \{ h \ | \  h \text{ is a Riemannian metric on } M \text{ and } h \text{ is } \varepsilon-\text{quasi-isospectral to } g\}.$$
Let $h_1, h_2 \in QI_\varepsilon(g)$, then we have that:

\begin{enumerate}
\item If $n$ is odd and $h_1$ and $h_2$ are $(\iota, \varepsilon)$-quasi-isospectral, then they are isospectral.
\item If $ n $ is even then
\begin{align*}
a_j(h_1) = a_j(h_2),& \,\text{ for }\, j < 1+n/2, \text{ and}\\
|a_{j}(h_1) - a_{j}(h_2)| \leq 2\varepsilon, & \,\text{ for }\, j=1+n/2.
\end{align*}
In addition, if $h_1$ and $h_2$ are $(\iota,\varepsilon)$-quasi-isospectral to $g$, then there is a constant $K=K(\lambda_\iota, \varepsilon)>0$, depending on $\lambda_\iota = \lambda_\iota(g)$ and $\varepsilon$,  such that
\begin{equation}|a_j(h_1) - a_j(h_2)| \leq 2\varepsilon, \ \text{ for } \  j > K. \label{eq:etail1} \end{equation}
\end{enumerate}
\label{thm:qiRmnb-hi}
\end{theorem}

Note that, in the statement of this theorem, we allow \(h_1\) and \(h_2\) to differ from \(g\) at different eigenvalues. Accordingly, the \((\iota,\varepsilon)\)-quasi-isospectral condition is not built into the definition of \(QI_\varepsilon(g)\).

\begin{proof}
Let $h_1, h_2 \in QI_\varepsilon(g)$, then for some $n_1, n_2 \geq 1$ we have that

\begin{eqnarray*}
\Spec(\Delta_{h_1}) & = & \{\lambda_m(g) + s_{1}\delta_{m,n_1} \}, \  \text{ for some } \  s_1 \in [-\varepsilon, \varepsilon], \\
\Spec(\Delta_{h_2}) & = &  \{\lambda_m(g) + s_{2}\delta_{m,n_2} \}, \  \text{ for some } \ s_2 \in [-\varepsilon, \varepsilon].
\end{eqnarray*}

From the spectral theory of differential operators on closed manifolds, we know that the following operators
$$e^{-t\Delta_{h_1}}, \quad e^{-t\Delta_{h_2}}, \quad e^{-t\Delta_{h_1}}-e^{-t\Delta_{h_2}} $$
are trace class and the traces have asymptotic expansions for small values of $t$, $t>0$, of the form given in equation~\ref{eq:aehtLnb} above.
Recall that the trace of a sum of trace-class operators is the sum of the traces. Moreover, we have that if all asymptotic expansions exist, we can sum the terms of the expansions, see Appendix~\ref{AsympExp}. Therefore, we obtain that for any $t >0$,
\begin{align*}
 \Tr\left(e^{-t \Delta_{h_1}}\right) - \Tr\left(e^{-t \Delta_{h_2}}\right) &= \sum_{\ell =0}^{\infty} e^{-t \lambda_{\ell}(h_1)} - \sum_{\ell =0}^{\infty} e^{- t \lambda_{\ell}(h_2)}\\
 &= e^{-t \lambda_{n_1}(h_1)} - e^{-t \lambda_{n_2}(h_2)} + \sum_{\ell \neq n_1}^{\infty} e^{-t \lambda_{\ell}(h_1)} - \sum_{\ell \neq n_2}^{\infty} e^{- t \lambda_{\ell}(h_2)}\\
& = e^{-t \lambda_{n_1}(h_1)} - e^{-t \lambda_{n_2}(h_2)} + e^{-t \lambda_{n_2}(h_1)} - e^{- t \lambda_{n_1}(h_2)}\\
&= e^{-t (\lambda_{n_1}(g) + s_1)} - e^{-t (\lambda_{n_2}(g) + s_2)} + e^{-t \lambda_{n_2}(g)} - e^{- t \lambda_{n_1}(g)}\\
&= e^{- t \lambda_{n_1}(g)}\left( e^{-t s_1} - 1\right) + e^{- t \lambda_{n_2}(g)}\left( 1- e^{-t s_2} \right).
\end{align*}
%%%%%%%%%%%%%%%%%%%%%%%%%%%%%%%%%%%%%%%%%%%%%%%
Let us consider two cases here, $s_1>0$ and $s_1<0$.
\begin{itemize}
\item {Case $s_1>0$:} We can easily see that $|e^{-t s_1} - 1| \leq t s_1 \leq \varepsilon t $.
\item {Case $s_1<0$:} In this case $t s_1 < 0$. We can proceed in two different forms. In the first form, we go a step backwards and recall that $\lambda_{n_1}(g) -\varepsilon >0$. Then
\begin{eqnarray*} e^{- t (\lambda_{n_1}(g) + s_1)} - e^{- t \lambda_{n_1}(g)}  &=&   e^{- t (\lambda_{n_1}(g) + s_1)}(1 - e^{t s_1}) \leq e^{- t (\lambda_{n_1}(g) + s_1)} |s_1| t \\  & \leq & e^{- t (\lambda_{n_1}(g) -\varepsilon)}  \varepsilon t.
\end{eqnarray*}
The second way in which we could proceed is the following: we assume that $|ts_1| \leq 1$. It is easy to prove that for  $x\in [0,1]$, $x \leq e^x -1 \leq 2x$. Therefore, we have that
$$e^{- t \lambda_{n_1}(g)}\left(e^{-t s_1} - 1\right) \leq  2 e^{- t \lambda_{n_1}(g)} |s_1| t \leq 2 \varepsilon \  t \ e^{- t \lambda_{n_1}(g)}.$$
Either of these options is good for us because we are interested in obtaining a linear estimate in $\varepsilon$ and a $O(t)$ as  $t\to 0^+$.  In what follows, we use the first form.
\end{itemize}

%%%%%%%%%%%%%%%%%%%%%%%%%%%%%%%%%%%%%%%%%%%%%%%%%
Proceeding in analogous way for $\lambda_{n_2}$ and $s_2$, we obtain that
\begin{eqnarray*}
\left| \Tr\left(e^{-t \Delta_{h_1}}\right) - \Tr\left(e^{-t \Delta_{h_2}}\right)\right|
&\leq&   e^{- t (\lambda_{n_1}(g) -\varepsilon)} t \varepsilon + e^{- t (\lambda_{n_2}(g)-\varepsilon)}  t \varepsilon \\
&\leq& \max \left\{ e^{- t (\lambda_{n_1}(g)-\varepsilon)}, e^{- t (\lambda_{n_2}(g)-\varepsilon)}\right\}  2 t \varepsilon \leq 2  t  \varepsilon.
\end{eqnarray*}

In this way, for small values of $t$, we have that
$$\Tr\left(e^{-t \Delta_{h_1}}\right) - \Tr\left(e^{-t \Delta_{h_2}}\right)  =  t^{-n/2} \sum_{j=0}^{\infty} b_j(h_1,h_2) t^j  = O(t).$$
The previous equations imply that $b_j(h_1,h_2)=0$ for  $j < 1+\frac{n}{2}$. Hence,
$$ a_j(h_1) = a_j(h_2) \quad \text{for }\quad  j < 1+\frac{n}{2}.$$
Now, for the sake of simplicity, let us restrict to the case when $n_1 = n_2 = \iota$, then
$$\lambda_{\iota}(h_1) = \lambda_{\iota}(g) + s_1 \ \text{ and } \lambda_{\iota}(h_2) = \lambda_{\iota}(g) + s_2,$$
and denote $\lambda_{\iota}(g) =\lambda_{\iota}$. Therefore,
$$\Tr\left(e^{-t \Delta_{h_1}}\right) - \Tr\left(e^{-t \Delta_{h_2}}\right) = e^{-t (\lambda_{\iota} + s_1)} - e^{-t (\lambda_{\iota} + s_2)}$$
The right-hand side is analytic and has a Taylor expansion at $t=0$. So, for $t \to 0^+$ we have
\begin{equation}
t^{-n/2} \sum_{j=0}^{\infty} b_j(h_1,h_2) t^j = \sum_{k=1}^{\infty} (-1)^k\frac{(\lambda_{\iota} + s_1)^k - (\lambda_{\iota} + s_2)^k}{k!}  t^k.
\label{eq:earhtEtedef}
\end{equation}
Since we have established above that the coefficients corresponding to $j< 1+\frac{n}{2}$ vanish, this equation becomes
\begin{equation}
\sum_{j = \lceil1+\frac{n}{2}\rceil}^{\infty} b_j(h_1,h_2) t^{j-n/2} = \sum_{k=1}^{\infty} (-1)^k\frac{(\lambda_{\iota} + s_1)^k - (\lambda_{\iota} + s_2)^k}{k!}  t^k,
\label{eq:earhtEtedef-pp}
\end{equation}
where $\lceil \ \cdot \ \rceil$ denotes the ceiling function.

If $n$ is odd, the left-hand side has only half-integer powers of $t$, while the right-hand side has only integer powers of $t$. This implies that all the coefficients vanish, i.e.,
\begin{align*}
b_j(h_1,h_2) &= 0, \ \text{ for } \ j\geq  \left\lceil1+\frac{n}{2}\right\rceil \\
(\lambda_{\iota} + s_1)^k - (\lambda_{\iota} + s_2)^k &= 0, \ \text{ for } k\geq 1.
\end{align*}
The equality of the heat invariants follows from the first equation. Since $\lambda_{\iota} >0$, we obtain $s_1 = s_2$. Hence, $h_1$ and $h_2$ are isospectral.

Now, let us consider an even-dimensional closed manifold, so $n$ is even. In this case, all powers of $t$ on both sides of equation (\ref{eq:earhtEtedef-pp}) are positive integers, and we have
\begin{equation}
a_j(h_1) - a_j(h_2) = (-1)^k\frac{(\lambda_{\iota} + s_1)^k - (\lambda_{\iota} + s_2)^k}{k!},
\end{equation}
for $k = j - \frac{n}{2} \geq 1$. Taking $k=1$, $j=1+ \frac{n}{2}$, and we obtain that
$$a_{1+ \frac{n}{2}}(h_1) - a_{1+ \frac{n}{2}}(h_2) = -(s_1 - s_2), $$
from which equation (\ref{eq:etail1}) follows.

Going back to the general case, we want to apply the Mean Value Theorem to the function $f(x) = (\lambda_{\iota} + x)^k$, w.l.o.g., assume that $s_2 < s_1$. Recall that $s_1, s_2 \in [-\varepsilon, \varepsilon]$. Therefore, there exist $w\in ]s_2, s_1[$ for which
\[ (\lambda_{\iota} + s_1)^k - (\lambda_{\iota} + s_2)^k = k (\lambda_{\iota} + w)^{k-1} (s_1 - s_2).\]
In this way
\begin{align}
a_j(h_1) - a_j(h_2) & = (-1)^{k} \frac{(\lambda_{\iota} + w)^{k-1}}{(k-1)!} (s_1 - s_2), \notag\\
| a_j(h_1) - a_j(h_2) | & \leq \frac{(\lambda_{\iota} + \varepsilon)^{k-1}}{(k-1)!} 2 \varepsilon, \label{eq:ubrhi-qimgc}
\end{align}
assuming that $\lambda_{\iota} - \varepsilon > 0$. The same estimate on the absolute value can be achieved if $s_1 < s_2$. Recall now that for any $a>0$
$$\lim_{k\to \infty} \frac{a^k}{k!} =0.$$
Using this and the fact that $\lambda_\iota$ and $\varepsilon$ are fixed, we obtain
$$\frac{(\lambda_{\iota} + \varepsilon)^{k}}{k!} \to 0 \  \text{ as }  \  k\to \infty.$$

Furthermore, consider the function $\varphi(x) = (\lambda_{\iota} + \varepsilon)^x / \Gamma(x+1)$, for which $\varphi(k) = \frac{(\lambda_{\iota} + \varepsilon)^{k}}{k!}$. A simple calculation using the properties of the digamma function $\psi=\Gamma'/\Gamma$,
shows that $\varphi$ attains a maximum at the value $x_0$ that satisfies
$$\log(\lambda_{\iota} + \varepsilon) = \psi(x_0 + 1),$$
 and for $x > x_0$, $\varphi(x)$ decreases.
Therefore, there is a $K>0$ such that for $k\geq K$
 $$\frac{(\lambda_{\iota} + \varepsilon)^{k-1}}{(k-1)!} \leq 1.$$
Notice that the dependence of $K$ on $\varepsilon$ does not affect the linearity in $\varepsilon$ of the estimate given in equation (\ref{eq:etail1}).
\end{proof}
We would like to point out that equation (\ref{eq:ubrhi-qimgc}) yields a uniform estimate for metrics in $QI_{\iota,\varepsilon}(g)$ which depends on $\lambda_\iota$, $\varepsilon$, and $k$. Furthermore, the estimate becomes uniform in $k$ and $j$ if one considers the maximum of the function $\varphi(x)$ on $]0,\infty[$. However, such a value may become very large if $\lambda_\iota$ is big.

Now, let us consider $(N,h)$ a Riemannian manifold with boundary. Let $\Delta_{h,D}$ be the Dirichlet Laplacian associated to the metric $h$ with domain
$H^1_0(M,h) \cap H^2(M,h)$. In this case, the asymptotic expansion of the heat trace as $t\to 0^+$ has the form:
\begin{equation}
\Tr\left(e^{-t\Delta_{h,D}}\right) = t^{-n/2} \sum_{j=0}^{\infty} a_j(h) t^{j/2} \label{eq:aehtLwb}
\end{equation}

\begin{proposition}
Let $N$ be a manifold with boundary of dimension $n$. Let $g$ be a Riemannian metric on it and let $\Delta_{g,D}$ be the associated Dirichlet Laplace operator. Set
$$QI_\varepsilon(g) = \{ h \ | \  h \text{ is a Riemannian metric on } N \text{ and } \Delta_{h,D} \text{ is } \varepsilon-\text{quasi-isospectral to } \Delta_{g,D}\}.$$
If $h_1, h_2 \in QI_\varepsilon(g)$.
Then,
$a_j(h_1) = a_j(h_2),$ if any of the following cases hold:
\begin{enumerate}
\item $j \leq n+1$.
\item If $n$ is even and $j$ is odd.
\item If $n$ is odd and $j$ is even.
\end{enumerate}
In addition,
$$|a_j(h_1) - a_j(h_2)| \leq 2\varepsilon, \ \text{ for } \ j = n+2,$$
and if $h_1$ and $h_2$ are $(\iota,\varepsilon)$-quasi-isospectral, then there exists a constant $K=K(\lambda_\iota, \varepsilon)>0$, such that equation (\ref{eq:etail1}) for $j > K$ whenever $j$ is not included in any of the preceding cases.
\label{eq:prop73}
\end{proposition}

\begin{proof}
The proof follows the same lines as the proof of Theorem ~\ref{thm:qiRmnb-hi}.

The main difference here is the boundary and the additional terms in the asymptotic expansion of the heat trace. In this case, instead of equation (\ref{eq:earhtEtedef}) we have
$$t^{-n/2} \sum_{j=0}^{\infty} b_j(h_1,h_2) t^{j/2} = \sum_{k=1}^{\infty} (-1)^k\frac{(\lambda_{\iota} + s_1)^k - (\lambda_{\iota} + s_2)^k}{k!}  t^k.$$
Therefore, all the coefficients with half-integer powers of $t$ vanish. If $n$ is even $n=2r$, for some $r\in \N$, and if $j$ is odd, $j=2s+1$ for some $s\in \N$, and
$$\frac j 2 - \frac n 2 = s - r +\frac 1 2 \notin \N.$$
When $n$ is odd, this will happen for $j$ even.

\end{proof}

Another important quantity associated with the spectrum of an operator is its regularised determinant. This is a number obtained from the spectrum via a process called zeta regularization. One starts by considering the spectral zeta function
$$\zeta_{\Delta_g} (s) = \sum_{\lambda_i > 0} \lambda_i^{-s}(g).$$
It converges uniformly and absolutely on compact subsets of $\Re(s)>\frac{n}{2}$. Using some properties of the heat trace, one can prove that there is a meromorphic extension of the spectral zeta function to $\C$ that is regular at $s=0$. Then the regularised determinant can be defined as
$$\det ( \Delta_g) = \exp(-\zeta'_{\Delta_g} (0)),$$
see \cite{rays}. Isospectral operators have the same determinant, and this is used in the proof of compactness of isospectral metrics on closed surfaces by Osgood et al. in \cite{OPS1}.
\begin{remark}
Let $g$ and $h$  be $(k,\varepsilon)$-quasi-isospectral Riemannian metrics defined on the same manifold. Then
$$\dfrac{\det( \Delta_g)}{\det(\Delta_h)} = \dfrac{\lambda_k(g)}{\lambda_k(h)}.$$
\end{remark}
An easy way to see this statement is by considering relative determinants. They were introduced by Werner M\"uller in \cite{Mueller:RD} using relative spectral zeta functions.
Let $g$ and $h$ be as above. Since each of the series defining the corresponding spectral zeta function converges absolutely and uniformly on compact subsets of $\Re(s)>\frac{n}{2}$, we have that
$$\zeta_{\Delta_g, \Delta_h} (s) = \sum_{i=1}^\infty \left( \lambda_i^{-s}(g)- \lambda_i^{-s}(h) \right) = \lambda_k(g)^{-s} - \lambda_k(h)^{-s}.$$
Notice that the series above actually converges for all $s\in \C$. Thus,
$$-\log\left(\det ( \Delta_g, \Delta_h) \right) = \zeta'_{\Delta_g, \Delta_h} (0) = -\log(\lambda_k(g)) + \log(\lambda_k(h)), $$
hence
$$\log(\det( \Delta_g, \Delta_h)) = \log(\lambda_k(g)) - \log(\lambda_k(h)), \quad \quad
\det( \Delta_g, \Delta_h) = \dfrac{\lambda_k(g)}{\lambda_k(h)}. $$
On the other hand, since the regularised determinant is defined for each of the operators, we also have
$$\det( \Delta_g, \Delta_h)  = \dfrac{\det( \Delta_g) }{\det(\Delta_h) }.$$
Putting this together, we obtain the claim in the remark.

\section{Quasi-isospectral potentials on Manifolds}
\label{sec:Results}
The results in this section are original, except where explicitly stated.
Let $(M,g)$ be a closed connected Riemannian manifold of dimension $n$. Let $V\in \cC^{\infty}(M)$ and $\Delta + V$ the corresponding Schr\"odinger operator. Let us assume $V\geq 0$. Let $\Spec(\Delta+V) = (\nu_m)_{m\geq 0}$ denote the spectrum of $\Delta + V$, notice that $\nu_0 \geq 0$. This is the same setting as in Section~\ref{subs:cisp}.

\begin{proposition}
Let $(M,g)$ be a closed Riemannian manifold of dimension $n$ and $V$ be a nonzero element in $\cC^{\infty}(M)$. Denote by $ (\nu_m)_{m\geq 0}$ the spectrum of $\Delta+V$. Let $k \geq 1$ and $\varepsilon>0$ be fixed. Let $(V_\ell)_{\ell\in \N}$ be a sequence of non-negative smooth functions such that for each $\ell$, the Schr\"odinger operator $\Delta + V_\ell$ is $(k,\varepsilon)$-quasi-isospectral to $\Delta+V$ with only the $k$-th eigenvalue being different, i.e. if $m\neq k$, $\lambda_m(V_\ell) = \nu_m$ and
$$\lambda_k(V_\ell) \in ]\nu_k- \varepsilon, \nu_k+\varepsilon[$$
for all $\ell \in \N$. Moreover, we require
$$\nu_{k-1} < \nu_k -\varepsilon \quad \text{ and } \quad \nu_k + \varepsilon < \nu_{k+1}.$$
Then, we have that:

\begin{enumerate}
\item If $n$ is odd, $V_\ell$ is isospectral to $V$, for each $\ell$.
\item If $ n $ is even, then
\begin{align*}
a_j(V) = a_j(V_\ell),& \,\text{ for }\, j < 1+n/2, \\
|a_{j}(V) - a_{j}(V_\ell)| \leq \varepsilon, & \,\text{ for }\, j=1+n/2,
\end{align*}
and there exists a constant $K=K(\nu_k, \varepsilon)>0$, depending on $\nu_k = \lambda_k(V)$ and $\varepsilon$, such that
$$|a_j(V) - a_j(V_\ell)| \leq \varepsilon, \ \text{ for } \ j >K.$$
\end{enumerate}
\label{th-caract_heat_coeff}
\end{proposition}

\begin{proof}
Then, for each $\ell$ and each $t>0$ we have that the heat operators
$$e^{-t(\Delta+V_\ell)} , \quad e^{-t(\Delta+V)} $$
are trace class and the corresponding traces have asymptotic expansions for small values of $t$ as given in equation (\ref{eq:aeHtSoSt}). Moreover, using Duhamel's principle one can also prove that $e^{-t(\Delta+V_\ell)} - e^{-t(\Delta + V)}$ is trace class, and it also has an asymptotic expansion as $t\to 0^+$ of the form given in equation (\ref{eq:ae-rht-st-sch}) since the same proof works in this setting. In the same way as in the proof of Theorem~\ref{thm:qiRmnb-hi}, we have that

\begin{align*}
 \Tr\left(e^{-t(\Delta+V_\ell)}\right) - \Tr\left(e^{-t(\Delta + V)}\right) &= \sum_{j =0}^{\infty} e^{-t \lambda_{j}(V_\ell)} - \sum_{j =0}^{\infty} e^{- t \nu_{j}}\\
 &= e^{-t \lambda_{k}(V_n) } - e^{-t \nu_{k}} + \sum_{j \neq k}^{\infty} e^{- t \lambda_j} - \sum_{j \neq k}^{\infty} e^{-t \nu_j} \\ &= e^{-t \lambda_{k}(V_n) } - e^{-t \nu_{k}}.
\end{align*}
Since $\lambda_k \in ]\nu_k- \varepsilon, \nu_k+\varepsilon[$, there is $s\in ]-\varepsilon, \varepsilon [$ such that $\lambda_k = \nu_k + s$. Using this and the same argument as in the proof of Theorem~\ref{thm:qiRmnb-hi}, we have that
\[\left|\Tr\left(e^{-t(\Delta+V_\ell)}\right) - \Tr\left(e^{-t(\Delta + V)}\right)\right|
= |e^{-t(\nu_k +s)} -  e^{-t \nu_{k}}| = e^{-t \nu_{k}}|e^{-ts}-1| \leq e^{-t (\nu_{k}-\varepsilon)} \varepsilon t
\]
and we obtain that, as $t\to 0$,
\[\Tr\left(e^{-t(\Delta+V_\ell)} - e^{-t(\Delta + V)}\right) =O(t).\]
On the other hand,
\begin{equation*}
\Tr\left(e^{-t(\Delta+V_\ell)} - e^{-t(\Delta + V)}\right) = \Tr\left(e^{-t(\Delta+V_\ell)}\right) - \Tr\left(e^{-t(\Delta + V)}\right),
\end{equation*}
and the asymptotic expansions of the l.h.s. and of the r.h.s. as $t\to 0$ are equal.
This implies, in particular, that
$$a_j(V_\ell) = a_j(V)$$
for all $j\geq 0$ if $n$ is odd, or for $j < 1+\frac{n}{2}$ if $n$ is even.

The rest of the proof follows the same lines as the proof of Theorem~\ref{thm:qiRmnb-hi}.
\end{proof}

In the same way as above, for the regularised determinants of the operators, we have that the relative determinant satisfies
$$\det( \Delta+V_\ell, \Delta+V)  = \dfrac{\det( \Delta+V_\ell) }{\det(\Delta+V)},$$
and we obtain the following expression:
$$\dfrac{\det( \Delta+V_\ell)}{\det(\Delta+V)} = \dfrac{\lambda_k(V_\ell)}{\lambda_k(V)} .$$

In Section~\ref{subs:cisp}, we discussed the results of Br\"uning \cite{Bruning1984OnTC} and Donnelly \cite{Donnelly2004CompactnessOI} regarding the compactness of isospectral potentials. In what follows, we show that analogous results hold in the case of quasi-isospectral potentials. Let $V_0 \in \mathcal{C}^\infty(M)$ and let $\mathcal{QI}_{k,\varepsilon}(V_0)$ denote the set of $(k,\varepsilon)$-quasi-isospectral potentials to $V_0$.

Thanks to Proposition~\ref{th-caract_heat_coeff}, we can follow the same argument used in \cite{Donnelly2004CompactnessOI} to conclude the compactness of $\mathcal{QI}_{k,\varepsilon}(V_0)$.

As mentioned before, the equality of the heat invariants for isospectral potentials implies a uniform bound in the norm of potentials in the Sobolev space $W_{n_0,2}(M)$, where
\[
n_0 := \min \left\{ s \in \mathbb{Z}_{\geq0} \,\middle|\, s > \frac{n}{2} - 2 \right\}.
\]

In each of the following cases, the proof follows directly from the proof of the corresponding result in \cite{Donnelly2004CompactnessOI} and the description of the relations between the heat invariants given in Proposition~\ref{th-caract_heat_coeff}. Note that in odd dimensions, we have that quasi-isospectrality implies isospectrality. Therefore, it remains only to establish the estimates in even dimensions.

\begin{proposition} \label{3-compactness}
    If $\dim(M) \leq 3$ and $V_0 \in C^\infty(M)$, then $\mathcal{QI}_{k,\varepsilon}(V_0)$ is compact in the $C^\infty(M)$ topology.
\end{proposition}

In this case, the first three invariants are required in order to obtain a uniform estimate for $\| V \| _{0,2} $, see \cite[Corollary 4.2]{Donnelly2004CompactnessOI}. In dimension $n=2$, the first two heat invariants coincide for quasi-isospectral potentials, and the difference between the third heat invariants is uniformly bounded by $\varepsilon$. Therefore, the inequalities still hold, and one can prove the necessary estimate.

Now, denoting $\mathcal{QI}^{+}_{k,\epsilon}(V_0)$ as the set of non-negative $(k,\epsilon)$-quasi-isospectral potentials to $V_0 \in C^\infty(M)$. In this space, we get a higher-dimensional result:
\begin{proposition}
    If $\dim(M) \leq 9$ and $V \in C^\infty(M)$, then $\mathcal{QI}^{+}_{k,\epsilon}(V_0)$ is compact in the $C^\infty(M)$ topology.
\end{proposition}
Following the proof of Proposition~5.1 in \cite{Donnelly2004CompactnessOI}, the first four heat invariants, up to $a_3$, are used to obtain a bound in $W_{1,2}(M)$ for dimensions $4$ and $5$. For dimensions $6$ and $7$, the fifth invariant $a_4$ is used to obtain a bound in $W_{2,2}(M)$. Finally, in dimensions $8$ and $9$, the equality of the sixth heat invariants, $a_5$, is used to prove a uniform bound in $W_{3,2}(M)$. As in the case of dimension~$2$, we do not have a proof that the necessary heat invariants coincide in each case, but they are uniformly bounded by $\varepsilon$, and this suffices to obtain the result.

Let us now consider dimension $n=1$ and Sturm-Liouville operators on the interval $[0,1]$ with Dirichlet boundary conditions.  Before we state our result, let us recall the following asymptotic expansion, which we write as a corollary because it follows straightforwardly from the results in \cite[Thm. 3.4]{Gilkey1979RecursionRA} and \cite{BanuelosBarreto}.

\begin{corollary}[from \cite{Gilkey1979RecursionRA} and \cite{BanuelosBarreto}]
Let $q\in \cC^{\infty}$, then the heat trace of $\cT_q$ (with Dirichlet boundary conditions) has the following asymptotic expansion as $t\to 0^+$
\begin{multline*}
\Tr \left(e^{-t \cT_q} \right) = \frac{1}{2\sqrt{\pi t}} - \frac{1}{2} + \frac{t^{1/2}}{\sqrt{4\pi}} \int_0^1 q(x)\; dx  + \frac{t}{4} (q(0)+q(1)) \\
+ \frac{t^{3/2}}{\sqrt{4\pi}} \left( \frac{1}{3} (q'(0)+q'(1)) + \int_0^1 q(x)^2 \; dx \right)  + O(t^2).
\end{multline*}
\end{corollary}

In this setting, we have the following Corollary:

\begin{corollary} Let $q\in {\cC^{\infty}[0,1]}$ and let $\cT_q$ the corresponding Sturm-Liouville operator with Dirichlet boundary conditions. Let $p$ be another potential such that $\cT_q$ and $\cT_p$ are $\varepsilon$-quasi-isospectral, for some $\varepsilon >0$. Then for $0\leq j\leq 2$ or $j$ even, we have that $a_j(\cT_p) = a_j(\cT_q)$. In particular
\begin{align*}
\int_0^1 q(x) dx &= \int_0^1 p(x) dx,\\
\frac{1}{3} (q'(0)+q'(1)) + \int_0^1 q(x)^2 \; dx &= \frac{1}{3} (p'(0)+p'(1)) + \int_0^1 p(x)^2 \; dx,\\
|q(0)+q(1) - (&p(0) + p(1))| \leq 2\varepsilon.
\end{align*}
\label{QI-sa}
\end{corollary}
The corollary follows from the computation done in the proof of Proposition~\ref{eq:prop73}. We show that for any $t > 0$
$$\Tr\left(e^{-t \cT_p}\right) - \Tr\left(e^{-t\cT_q}\right) = e^{-t \lambda_{k}(p) } - e^{-t \lambda_{k}(q)}.$$
Therefore, the difference of their asymptotic expansions is also $O(t)$ as $t\to 0^+$. However, as we pointed out above, since we are in dimension $n=1$, that is odd, and the potentials are smooth, we have that the heat invariants are the same, $a_j(p) = a_j(q)$ for $j \leq 2$ or $j$ even.

Finally, observe that the arguments presented in this section extend to the case in which the spectra of two operators differ in not just one, but in finitely many eigenvalues.

\section{Appendix}

\subsection{Asymptotic expansions}
\label{AsympExp}
The main argument of this paper relies on the fact that we can manipulate asymptotic expansions. This involves exchanging several limit processes. As any person who has taken a course in mathematical analysis knows, sometimes it is clear that one can do it, however other times exchanging limits may lead to false statements. In any case, it is always necessary to make sure that one is allowed to perform these manipulations by carefully justifying each step.~In this appendix, we include the definition and some basic facts about asymptotic expansions. We follow the books by De Brujin \cite{Bruijn2014asymptotic} and by Miller \cite{miller2006applied}

\begin{definition} Let $D$ be a domain. A sequence  of functions $(\varphi_j(t))_{j\in \N}$  on $D$ is called an asymptotic sequence (scale) as $t \to t_0$ from $D$ if for all $n>m$,
$$\varphi_{n}(t) = o(\varphi_{m}(t)) \ \text{as} \  t\to t_{0}, \ t\in D.$$
In particular,
$\varphi_{n+1}(t) = o(\varphi_{n}(t))$ as $t\to t_{0}$, $t\in D$.
\end{definition}

\begin{example} Let $\alpha >0$, then $(\varphi_j(t) = t^{-\alpha + j})_{j\in \N}$ forms an asymptotic sequence as $t\to 0^+$ for $t \in ]0,\infty[$.
\end{example}

\begin{definition}
We say that a function $f$ has an asymptotic expansion with respect to the asymptotic sequence $(\varphi_j(t))_{j\in \N}$ as $t\to t_0$, and we write
$$f(t) \sim \sum_{j=0}^{\infty} a_j \varphi_j(t)$$
if for every $N\in \N$, we have that
$$f(t)-\sum_{j=0}^N a_j \varphi_j(t) = o(\varphi_N (t)) \text{ as } t\to t_0. $$
\end{definition}
Let us list some facts about asymptotic series as $t\to t_0$ from $D$:
\begin{enumerate}
\item An asymptotic series is not necessarily a convergent series.
\item If $f(t) \sim \sum_{j=o}^{\infty} c_j \varphi_j(t)$, the coefficients are unique.
\item A series $\sum_{j=o}^{\infty} c_j \varphi_j(t)$ may approximate different functions, for example,
$$e^{-t} \sim 0 \cdot 1 +  0\cdot t^{-1} + 0 \cdot t^{-2} + \dots,$$ as $t\to \infty$.
\end{enumerate}
The main property of asymptotic expansions for us is the following:
If
\begin{eqnarray*}
f(t) &\sim& \sum c_j \varphi_j(t),\\
g(t) &\sim& \sum b_j \varphi_j(t), \text{ and }\\
h(t):= (f+g)(t) &\sim& \sum a_j \varphi_j(t),
\end{eqnarray*}
as $t\to t_0$, then for all $j\in \N$, we have that
$$a_j = c_j + b_j.$$

\begin{remark}
Comparison with Taylor series. A Taylor series is a particular case of a power series and has a radius of convergence. Inside the interior of its disk of convergence, the Taylor series converges to the function. Note that the radius of convergence may be zero. Every power series defines an asymptotic expansion, but not every asymptotic expansion arises from a power series.
\end{remark}

\subsection{Proof of Darboux Lemma}

In this section, we include a proof of Darboux's Lemma as stated in Lemma~\ref{lemma:Darboux}.

\begin{proof}(Darboux lemma)
    The following argument is known as the commutation method \cite[Deift]{deift}. The trick is to decompose the operator
$$H=-\dfrac{d^2}{dx^2}+q-\mu$$
as a multiplication of other operators, and to compare it with its commutation.
Let
$$T=g\left(\dfrac{d}{dx}\right)\frac{1}{g} \text{ and }S=-\frac{1}{g}\left(\dfrac{d}{dx}\right)g.$$
Here $S$ is the formal adjoint (on $\cL^2(]0,1[)$) of $T$: Indeed, let $\phi, \psi \in \cL^2(]0,1[)$
\begin{align*}
    \langle T(\phi),\psi\rangle&=\int_0^1 T(\phi(x))\psi(x)\; dx=\int_0^1 \phi(x)\left[-\frac{1}{g(x)}\dfrac{d}{dx}(g(x)\psi(x))\right]\; dx=\langle \phi,S(\psi)\rangle.
\end{align*}
Here the integration by parts was used.
On the other hand, by calculating a composition between them:
\begin{align*}
   ST(\phi) &=-\frac{1}{g}\left(\dfrac{d}{dx}\left[g\dfrac{d\phi}{dx}-\dfrac{dg}{dx}\phi \right]\right)\\
    &=-\dfrac{d^2\phi}{dx^2}+(q-\mu)\phi = H(\phi),
\end{align*}
while in the other order:
\begin{align*}
    TS(\phi) &=g\left(\dfrac{d}{dx}\left[-\frac{1}{g^{2}}\dfrac{dg}{dx}\phi-\frac{1}{g}\dfrac{d\phi}{dx}\right]\right)\\
    &=-\dfrac{d^2\phi}{dx^2}-\mu \phi+\left[q-2\dfrac{d^2}{dx^2}\log(g)\right]\phi.
\end{align*}

According to the notation in the lemma, $f$, $g$ and $h$ satisfy
\[
    ST(f)=(\lambda-\mu)f, \quad
    ST(g)=0 \quad \text{and} \quad
    TS(h)=(\lambda-\mu)h.
\]
In particular,
$$Tf=g\left(\dfrac{d}{dx}\right)\frac{1}{g}f=-\frac{1}{g}\dfrac{dg}{dx}f+\dfrac{df}{dx}=\frac{1}{g}[g,f]=h$$
is a solution to $TST(f)=(\lambda-\mu)Tf$. On the other hand, $Tf$ is a non-trivial solution to $TST(f)=(\lambda-\mu)Tf$ if and only if $f$ is a non-trivial solution to $ST(f)=(\lambda-\mu)f$.
\end{proof}

%%% REFERENCES %%%

\end{document}